%% file: main.tex
\documentclass{article}

\usepackage{PRIMEarxiv}

\usepackage[utf8]{inputenc} % allow utf-8 input
\usepackage[T1]{fontenc}    % use 8-bit T1 fonts
\usepackage{hyperref}       % hyperlinks
\usepackage{url}            % simple URL typesetting
\usepackage{booktabs}       % professional-quality tables
\usepackage{amsfonts}       % blackboard math symbols
\usepackage{nicefrac}       % compact symbols for 1/2, etc.
\usepackage{microtype}      % microtypography
\usepackage{lipsum}
\usepackage{fancyhdr}       % header
\usepackage{graphicx}       % graphics
\graphicspath{{media/}}     % organize your images and other figures under media/ folder
\usepackage{xcolor}
\usepackage{mathtools}
\usepackage{bm}
\usepackage{subcaption}
\usepackage[square,numbers]{natbib}
\usepackage{soul}
\usepackage{tikz}
\usepackage{circuitikz}
\usepackage{xurl}

\usepackage{tablefootnote}

\usepackage{pifont}% http://ctan.org/pkg/pifont
\newcommand{\cmark}{\ding{51}}%
\usepackage{xcolor}

\usepackage[acronym]{glossaries}
\input{glossaries}
\glsdisablehyper
\input{definitions}

\usepackage{enumitem}
  
\title{A Frequency-Domain Differential Corrector for Quasi-Periodic Trajectory Design and Analysis
}

\author{
  Beom Park\\
  PhD Candidate \\
  School of Aeronautics and Astronautics \\
  Purdue University \\
  West Lafayette, IN, USA, 47906 \\
  \texttt{park1103@purdue.edu} \\
  \And
  Kathleen C. Howell \\
  Hsu Lo Distinguished Professor \\ 
  School of Aeronautics and Astronautics \\
  Purdue University \\
  West Lafayette, IN, USA, 47906 \\
  \texttt{howell@purdue.edu}
  \And
  Shaun Stewart \\
  Flight Dynamics Lead \\
  Intuitive Machines \\
  Houston, TX, USA, 77059\\
  \texttt{sstewart@intuitivemachines.com} \\
}

\begin{document}
\maketitle

\input{text/abstract.tex}
\input{text/intro.tex}

\input{text/dynamical_models}
\input{text/quasi_periodic_orbits_freq}
\input{text/single_shooter}
\input{text/multiple_shooter}
\input{text/conclusions}

\section*{Acknowledgments}
Valuable discussions with members from Multi-Body Dynamics Research Group are appreciated. Portions of this work were completed at Purdue University under Intuitive Machines LLC Award 21123283.

\bibliographystyle{custom}  
\bibliography{references}  

\clearpage

\input{text/partials}

\input{text/sensitivities}

\input{text/sensitivities_ms}

\end{document}

%% file: glossaries.tex
\newacronym{cr3bp}{CR3BP}{Circular Restricted Three-Body Problem}
\newacronym{er3bp}{ER3BP}{Elliptic Restricted Three-Body Problem}
\newacronym{bcr4bp}{BCR4BP}{Bi-Circular Restricted Four-Body Problem}
\newacronym{hr3bp}{HR3BP}{Hill Restricted Three-Body Problem}
\newacronym{hr4bp}{HR4BP}{Hill Restricted Four-Body Problem}
\newacronym{qbcp}{QBCP}{Quasi Bi-Circular Problem}
\newacronym{qhr4bp}{QHR4BP}{Quasi-Hill Restricted Four-Body Problem}
\newacronym{qhr4bps}{QHR4BPs}{Quasi-Hill Restricted Four-Body Problems}

\newacronym{iqhr4bp}{I-QHR4BP}{In-plane Quasi-Hill Restricted Four-Body Problem}
\newacronym{oqhr4bp}{O-QHR4BP}{Out-of-plane Quasi-Hill Restricted Four-Body Problem}

\newacronym{hfem}{HFEM}{Higher-Fidelity Ephemeris Model}

\newacronym{nrho}{NRHO}{Near Rectilinear Halo Orbit}

\newacronym{qpo}{QPO}{Quasi-Periodic Orbit}
\newacronym{qpos}{QPOs}{Quasi-Periodic Orbits}

\newacronym{myt}{MYT}{Multi-Year Trajectory}
\newacronym{myts}{MYTs}{Multi-Year Trajectories}

\newacronym{gmos}{GMOS}{G{\'o}mez Mondelo Olikara Scheeres}

\newacronym{po}{PO}{Periodic Orbit}
\newacronym{pos}{POs}{Periodic Orbits}

\newacronym{dft}{DFT}{Discrete Fourier Transform}
\newacronym{fft}{FFT}{Fast Fourier Transform}
\newacronym{cft}{CFT}{Continuous Fourier Transform}

\newacronym{dro}{DRO}{Distant Retrograde Orbit}
\newacronym{dros}{DROs}{Distant Retrograde Orbits}
\newacronym{stm}{STM}{State Transition Matrix}

\newacronym{qpts}{QPTs}{Quasi-Periodic Trajectories}

\newacronym{brf}{BRF}{Barycentric Rotating Frame}
\newacronym{eof}{EOF}{Earth Orbit Frame}
\newacronym{mci}{JIF}{J2000 Inertial Frame}

\newacronym{dam}{DADM}{Doubly-Averaged Dynamical Model}

\newacronym{elfos}{ELFOs}{Elliptical Lunar Frozen Orbits}
\newacronym{elfo}{ELFO}{Elliptical Lunar Frozen Orbit}

\newacronym{laskar}{L-NAFF}{Laskar-Numerical Analysis of Fundamental Frequency}
\newacronym{gomez}{GMS-C}{G{\'o}mez-Mondelo-Sim{\'o}-Collocation}

\newacronym{fdt}{FDT}{Frequency Domain Targeter}

\newacronym{fmp}{FMP}{Flow Map Parametrization}
\newacronym{fdc}{FDDC}{Frequency-Domain Differential Corrector}

%% file: definitions.tex
\renewcommand{\exp}[1]{e^{#1}}

\newcommand{\argmax}[1]{\underset{#1}{\text{arg max}}}

\newcommand{\traj}{\ensuremath{\gamma}}

\newcommand{\signalt}{\ensuremath{q}} % signal in the time domain
\newcommand{\signala}{\ensuremath{\tilde{q}}} % 
\newcommand{\signalaidxA}{\ensuremath{\tilde{q}}_{\idxA}} % approximated signal
\newcommand{\signalf}{\ensuremath{\mathcal{D}_{\signalt}}} % signal in the frequency domain
\newcommand{\signalcosf}{\ensuremath{\mathcal{D}_{c(\freqA{\idxA})}}} % signal in the frequency domain
\newcommand{\signalsinf}{\ensuremath{\mathcal{D}_{s(\freqA{\idxA})}}} % signal in the frequency domain

\newcommand{\cosf}{\ensuremath{\mathcal{C}_{\signalt}}}
\newcommand{\sinf}{\ensuremath{\mathcal{S}_{\signalt}}}
\newcommand{\csf}{\ensuremath{\mathcal{CS}_{\signalt}}}

\newcommand{\pvec}[1]{\vec{#1}\mkern2mu\vphantom{#1}}

\newcommand{\da}[1]{\ensuremath{\bar{\bar{#1}}}}

\definecolor{ao}{rgb}{0.0, 0.5, 0.0}
\definecolor{amber}{rgb}{1.0, 0.75, 0.0}
\definecolor{cyan}{rgb}{0.0, 1.0, 1.0}
\definecolor{magenta}{rgb}{1.0, 0.0, 1.0}
\definecolor{purple}{rgb}{0.5, 0.0, 0.5}

\newcommand{\yes}{\textcolor{ao}{\cmark}}

\def\thickhline{\noalign{\hrule height2.0pt}}

\newcommand{\lag}{\textcolor{amber}{\rotatebox[origin=c]{90}{\ding{220}}}}
\newcommand{\lead}{\textcolor{purple}{\rotatebox[origin=c]{270}{\ding{220}}}}

\newcommand{\coscosf}{\ensuremath{\mathcal{C}_{c(\freqA{\idxA})}}}
\newcommand{\cossinf}{\ensuremath{\mathcal{C}_{s(\freqA{\idxA})}}}
\newcommand{\sincosf}{\ensuremath{\mathcal{S}_{c(\freqA{\idxA})}}}
\newcommand{\sinsinf}{\ensuremath{\mathcal{S}_{s(\freqA{\idxA})}}}
\newcommand{\cscosf}{\ensuremath{\mathcal{CS}_{c(\freqA{\idxA})}}}
\newcommand{\cssinf}{\ensuremath{\mathcal{CS}_{s(\freqA{\idxA})}}}

\newcommand{\sampleSize}{\ensuremath{N}}

\newcommand{\nA}{\ensuremath{m}}

\newcommand{\ampA}[1]{\ensuremath{A_{#1}}}
\newcommand{\freqA}[1]{\ensuremath{\mathtt{f}_{#1}}}

\newcommand{\freqP}[1]{\omega_{\idxF}^{*#1}} % j-th peak
\newcommand{\freqPO}[1]{\omega_{\idxF\pm1}^{*#1}} % j-th peak and its nearby 2nd highest freq.
\newcommand{\freqPP}[1]{\omega_{\idxF+1}^{*#1}} % j-th peak and its nearby 2nd highest freq.
\newcommand{\freqPM}[1]{\omega_{\idxF-1}^{*#1}} % j-th peak and its nearby 2nd highest freq.
\newcommand{\freqPPM}[1]{\omega_{\idxF\pm1}^{*#1}} % j-th peak and its nearby 2nd highest freq.

\newcommand{\fund}[1]{\nu_{#1}}

\newcommand{\freqAArg}{\ensuremath{\mathtt{f}}}
\newcommand{\phaseA}[1]{\ensuremath{\theta_{#1}}}
\newcommand{\idxA}{\ensuremath{j}}

\newcommand{\idxF}{\ensuremath{k}} % freq index
\newcommand{\tidxT}{\tnd_{\idxT}}
\newcommand{\tspan}{\ensuremath{T}} % total span, but might have to change later

\newcommand{\idxT}{\ensuremath{i}} % idx for the transformation sum
\newcommand{\imag}{\ensuremath{\mathtt{i}}}
\newcommand{\han}{\ensuremath{h}}

\newcommand{\ex}{\ensuremath{n_e}}
\newcommand{\inter}{\ensuremath{n_i}}
\newcommand{\ttl}{\ensuremath{n}}

\newcommand{\freeVE}{\ensuremath{\vec{\xi}}}

\newcommand{\constraintV}[1]{\ensuremath{\vec{F}_{#1}}}
\newcommand{\constraintC}[1]{F_{#1}} % component
\newcommand{\constraintVL}{\constraintV{L, \idxA}}
\newcommand{\constraintVG}{\constraintV{G, \idxA}}

\newcommand{\tnd}{\ensuremath{t}}

\newcommand{\tndS}{\Delta \tnd} % sampled time step
\newcommand{\fnd}{\omega}
\newcommand{\fidxF}{\fnd_{\idxF}} 
\newcommand{\fidxFm}{\fnd_{\idxF-1}}
\newcommand{\fidxFp}{\fnd_{\idxF+1}}% sampled frequency
\newcommand{\fndS}{\Delta \fnd } % sampled freq step

\newcommand{\dft}[2]{\mathcal{D}_{#1}(#2)}
\newcommand{\genericFreq}{f}

\newcommand{\rt}[2]{\ensuremath{#1{\,:\,}#2}}

\newcommand{\x}[1]{\vec{\mathbb{X}}_{#1}}
\newcommand{\xp}[4]{\vec{\mathtt{X}}^{#2}(#3; #4)}
\newcommand{\xpp}{\vec{\mathtt{X}}}
\newcommand{\epoch}[1]{\mathtt{t}_{#1}}
\newcommand{\real}[1]{\mathbb{R}^{#1}}

\newcommand{\f}[1]{\vec{\mathbb{F}}_{#1}}

\newcommand{\stm}[4]{\bm{\phi}^{#2}(#3; #4)}

\newcommand{\nP}{\ensuremath{n_p}} % number of pps
\newcommand{\mb}{\begin{bmatrix}}
\newcommand{\me}{\end{bmatrix}}

%% file: text/abstract.tex
\begin{abstract}
This paper introduces the \acrfull{fdc}, a model-agnostic approach for constructing quasi-periodic orbits (QPOs) across a range of dynamical regimes. In contrast to existing methods that explicitly enforce an invariance condition in all frequency dimensions, the \acrshort{fdc} targets dominant spectral components identified through frequency-domain analysis. Leveraging frequency refinement strategies such as \acrfull{laskar} and \acrfull{gomez}, the method enables efficient and scalable generation of high-dimensional QPOs. The \acrshort{fdc} is demonstrated in both single- and multiple-shooting formulations. While the study focuses on the Earth–Moon system, the framework is broadly applicable to other celestial environments. Sample applications include Distant Retrograde Orbits (DROs), Elliptical Lunar Frozen Orbits (ELFOs), and Near Rectilinear Halo Orbits (NRHOs), illustrating constellation design and the recovery of analog solutions in higher-fidelity models. With its model-independent formulation and spectral targeting capabilities, \acrshort{fdc} offers a versatile tool for robust trajectory design and mission planning in complex dynamical systems.
\end{abstract}

%% file: text/intro.tex
\section{Introduction}

Periodic and quasi-periodic trajectories serve as a fundamental framework in astrodynamics. Closely linked to the integrability of a system, these structures illuminate the architecture of the phase space and provide a basis for understanding long-term orbital behaviors. For instance, the classical two-body problem, comprising a point-mass central body and a spacecraft, is a superintegrable system in the Liouville-Arnold sense \cite{arnol2013mathematical}, with five independent integrals of motion within a three degree-of-freedom Hamiltonian system. Its solutions emerge as the conic sections including the ellipses as the fundamental periodic structures, undoubtedly serving as basis for astrodynamics. When perturbations are introduced, often rendering the system non-integrable, new behaviors such as chaos may emerge. Nonetheless, many periodic and quasi-periodic structures persist. The \acrfull{cr3bp}, describing a spacecraft under the gravitational influence of two primary bodies, serves as a canonical example for such a non-integrable system. Typically represented within a rotating frame, it retains only a single integral of motion, the Jacobi constant \cite{szebehely2012theory}. Despite its non-integrability, the CR3BP admits a wealth of (quasi-)periodic orbits, including those near each primary and various resonant orbits, as well as families of orbits near the Lagrange points \cite{henon1997generating, doedel2007elemental}.

The (quasi-)periodic orbits play a critical role in both theoretical and applied contexts in astrodynamics. In particular, periodic orbits in the CR3BP exist in continuous one-parameter families \cite{henon1997generating}. Linearized dynamics around the periodic orbit allow bifurcation analysis \cite{brouckeER3BP1969} and aid in characterizing the invariant manifolds that govern the nearby phase space \cite{seydel2009practical, wiggins2003}. In the planar CR3BP, for example, the hyperbolic manifolds associated with Lyapunov orbits around the $L_1$ and $L_2$ Lagrange points are responsible for organizing global transport and giving rise to chaotic motion through homoclinic tangles \cite{koon2000dynamical, conley1968low, almanza2025hawaii}. Beyond their utility in illustrating the global phase space, specific (quasi-)periodic orbits are often well-suited for mission design, offering useful properties such as repeated ground tracks, reduced eclipse exposure, constant line-of-sight, or proximity to regions of interest, such as the lunar south pole. Notable mission-relevant examples include the Near Rectilinear Halo Orbit (NRHO) and the Distant Retrograde Orbit (DRO), underpinning key components of NASA’s Gateway architecture \cite{zimovan2022dynamical}.

Despite the utility of (quasi-)periodic orbits in the CR3BP, identifying suitable counterparts in higher-fidelity models remains a significant challenge. As noted by Jorba and Villanueva \cite{jorba1997persistence}, an 
$\inter$-dimensional quasi-periodic orbit (QPO) in the CR3BP generally evolves into an $(\ttl = \inter + \ex)$-dimensional QPO in a higher-fidelity model, with the additional $\ex$ frequencies inherited (external) from perturbations. For the CR3BP, $\inter \leq 3$, as a $3$ degree-of-freedom Hamiltonian system. However, the number of inherited frequencies $\ex$ may be significant. For example, within the Earth-Moon system, at least five distinct frequencies ($\ex \geq 5$) are typically required to capture key perturbations, including Earth-Moon distance variations and solar effects \cite{gomez2002solar}.

The current investigation focuses on developing numerical strategies to construct higher-dimensional quasi-periodic trajectories that occur frequenctly within more realistic dynamical environments. The increase in dimensionality poses challenges for existing numerical schemes. Baresi, Olikara, and Scheeres \cite{baresi2018fully} review several fully-numerical approaches, with an algorithm proposed by G{\'o}mez and Mondelo \cite{gomez2001dynamics2} and refined by Olikara and Scheeres \cite{olikara2012numerical} (GMOS) serving as the de facto standard due to its advantages in computational speed and accuracy. However, GMOS and similar strategies face the ``curse of dimensionality,'' where the cost of supplying high-dimensional QPOs scales poorly due to the need to invert large, dense matrices \cite{kumar2022rapid}. To date, this method has been verified primarily for $\ttl = \inter + \ex \leq 3$ \cite{mccarthy2022diss, iac2024, baresi_iac, brown2025hawaii}, and while extensions are possible, computational requirements quickly become intractable as $\ttl$ increases. An alternative approach proposed by Haro and Mondelo \cite{haro2021flow}, termed the \acrfull{fmp}, circumvents some of these limitations by replacing the dense matrix inversion with a quasi-Newton iteration. The method is semi-analytical and leverages quasi-Floquet theory \cite{jorba1992reducibility} for an alternative basis. Although FMP significantly reduces computational cost and has been demonstrated and extended to the \acrfull{er3bp} \cite{fernandez2024flow, kumar2022rapid} and the \acrfull{hr4bp} \cite{henry2024,Dennis2025Hawaii}, its implementation is more complex compared to the GMOS algorithm. The examples that appear in the literature discuss $\ttl \leq 3$. Whether the method scales effectively to larger values of $\ttl $ remains an open question. Both GMOS and FMP aim to explicitly construct the $\ttl $-dimensional QPO directly. Although the methods are formally generalizable, their reliance on model-specific derivations may pose a barrier to broader adoption and practical application.

Proposed in this paper is a complementary approach to deliver $\ttl$-dimensional QPOs, termed the \acrfull{fdc}. The strategy acknowledges that the $\ex$ inherited frequencies are dictated by external forcing and are not directly controllable. By leveraging frequency-domain information, specifically the correction of dominant spectral peaks, the strategy introduces a differentiable update process informed by previous work on frequency analysis \cite{laskar1999introduction, gomez2010collocation}. Thus, rather than explicitly enforcing the invariance condition in all $\ttl$ dimensions, the \acrshort{fdc} corrects for the intrinsic $\inter$ frequency components. The \acrshort{fdc} method is intended as a more computationally tractable strategy for constructing high-dimensional QPOs in complex dynamical environments. Moreover, the approach is model-agnostic: the same formulation is adaptable across a variety of dynamical models without requiring significant modification. These features make the \acrshort{fdc} strategy well-suited for generating baseline solutions for long-duration missions such as Gateway, or for maintaining constellations of spacecraft with constraints on relative phasing and geometry.

The remainder of the document is organized as follows. Section \ref{sec:dynamics} introduces representative dynamical models within the Earth-Moon system, although the proposed method is applicable to any other regimes as well. Foundational tools in frequency analysis are explored in Section \ref{sec:frequency_refinement} that motivate the corrections strategy. In Section \ref{sec:single}, the single-shooting formulation of the \acrshort{fdc} strategy is presented, along with demonstrations involving \acrfull{dros} and \acrfull{elfos}. The algorithm is extended in Section \ref{sec:multi} to a multiple-shooting formulation, with an application to a $\rt{9}{2}$ synodic resonant halo orbit. Concluding remarks are provided in Section \ref{sec:summary}.

%% file: text/dynamical_models.tex
\section{\label{sec:dynamics}Preliminaries: Frames, Dynamical Models, and \acrfull{qpos}}

The \acrfull{fdc} is designed to be model-agnostic, and its implementation is demonstrated using three representative dynamical models in the Earth–Moon system: the \acrfull{dam}, the \acrfull{cr3bp}, and a \acrfull{hfem}. These models span a hierarchy of fidelity levels and illustrate the scalability of the proposed approach. The following subsections detail the coordinate systems, frame transformations, and governing equations for each model.

\subsection{Frames and Transformations}

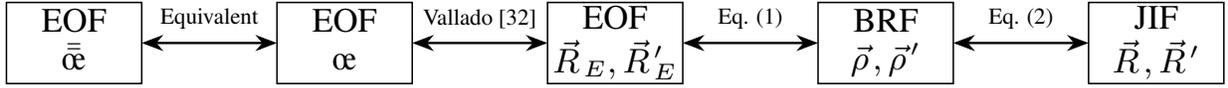
\begin{figure}[!ht]
\centering
\resizebox{1\textwidth}{!}{%
\begin{circuitikz}
\tikzstyle{every node}=[font=\small]

\draw  (1.25,9) rectangle  node {\footnotesize \begin{tabular}{c}
     \acrshort{eof}  \\
     $\bar{\bar{\text{\oe}}}$ \\
  \end{tabular}} (2.5,8.25);
\draw  (3.75,9) rectangle  node {\footnotesize  \begin{tabular}{c}
     \acrshort{eof}  \\
     \oe \\
  \end{tabular}} (5,8.25);
\draw  (6.25,9) rectangle  node {\footnotesize  \begin{tabular}{c}
     \acrshort{eof}  \\
     $\pvec{R}_E, \pvec{R}_E'$ \\
  \end{tabular}} (7.5,8.25);
\draw  (8.75,9) rectangle  node {\footnotesize  \begin{tabular}{c}
     \acrshort{brf}  \\
     $\pvec{\rho}, \pvec{\rho}'$ \\
  \end{tabular}} (10,8.25);
\draw  (11.25,9) rectangle  node {\footnotesize  \begin{tabular}{c}
     \acrshort{mci}  \\
     $\pvec{R}, \pvec{R}'$ \\
  \end{tabular}} (12.5,8.25);
\draw [line width=0.6pt, <->, >=Stealth] (2.5,8.625) -- (3.75,8.625);
\draw [line width=0.6pt, <->, >=Stealth] (5,8.625) -- (6.25,8.625);
\draw [line width=0.6pt, <->, >=Stealth] (7.5,8.625) -- (8.75,8.625);
\draw [line width=0.6pt, <->, >=Stealth] (10,8.625) -- (11.25,8.625);
\node [font=\tiny] at (3.125,8.85) {Equivalent};
\node [font=\tiny] at (5.625,8.85) {Vallado \cite{vallado2001fundamentals}};
\node [font=\tiny] at (8.125,8.85) {Eq. \eqref{eq:eof_brf}};
\node [font=\tiny] at (10.625,8.85) {Eq. \eqref{eq:mci_brf}};
\end{circuitikz}
 
}%
\caption{\label{fig:frame_transformation}Transformations across various state representations and frames}  

\end{figure}

To support analyses across these models, three reference frames are employed. Each frame plays a distinct role in state propagation, element extraction, and frame-to-frame transformation. These various frames facilitate analyses of various structures as follows:
\vspace{-\parskip}
\begin{enumerate}[itemsep=0mm,  label= \color{red}(Frame \arabic*), leftmargin=3mm, itemindent=20mm, labelsep=!, labelwidth=2mm, labelindent=20mm]
\item \ul{Barycentric Rotating Frame} (\acrshort{brf}): this frame adopts an orthogonal basis $\hat{x}-\hat{y}-\hat{z}$ defined as follows: (1) $\hat{x}$ is directed from the Earth to the Moon, (2) $\hat{z}$ coincides with the angular momentum of the Moon in its orbit with respect to Earth, and (3) $\hat{y}$ completes the dextral triad. The origin of the frame is located at the Earth-Moon barycenter. With the Earth-Moon distance normalized to be unity (regardless of the actual dimensional distance), this frame renders fixed locations for the Earth and Moon along the $\hat{x}$ axis at $-\mu$ and $(1-\mu)$. The mass ratio $\mu$ measures the ratio of the mass of the Moon with respect to the Earth-Moon system mass. The nondimensional (nd) spacecraft position vector is expressed as $\vec{\rho} = x\hat{x} + y\hat{y} + z\hat{z}$.
\item \ul{Moon-Centered Earth Orbit Frame} (\acrshort{eof}): the orthonormal basis $\hat{X}_E-\hat{Y}_E-\hat{Z}_E$ is constructed such that $\hat{Z}_E$ coincides with $\hat{z}$, and $\hat{X}_E-\hat{Y}_E$ are produced by rotating $\hat{x}-\hat{y}$ at a constant negative rate with respect to $\hat{Z}$. The origin of the frame is located at the center of the Moon. The nd spacecraft position vector within the \acrshort{eof} is then constructed as, \\
\begin{align}
    \label{eq:eof_brf} \vec{R}_E = \bm{C}_E (\vec{\rho}-\vec{\rho}_M) = X_E\hat{X}_E + Y_E\hat{Y}_E + Z_E\hat{Z}_E, \quad \bm{C}_E =\mb \cos (t ) & -\sin (t ) & 0 \\ \sin (t ) & \cos (t ) & 0 \\ 0 & 0 & 1 \me,
\end{align}
where $\rho_M = (1-\mu)\hat{x}$, the nd lunar position vector within the \acrshort{brf}. The nd time $t$ is constructed leveraging the characteristic quantities $l_*$ and $t_*$ defined as the average Earth-Moon distance and $t_* = \sqrt{l_*^3/(\tilde{\mu}_E + \tilde{\mu}_M)}$ where $\tilde{\mu}_{E, M}$ are the dimensional gravitational parameters for the Earth and Moon, respectively. Then, $t$ normalizes the dimensional time via $t_*$. It is assumed that, at $t = 0$, the Earth is on the positive $\hat{X}_E$-axis without losing generality. This frame is defined consistent with previous investigations \cite{broucke2003long, ely2005stable, folta2006lunar}. For the nd position and velocity, i.e., $\vec{R}_E, \pvec{R}_E' = d\vec{R}_E/dt$, the osculating orbital elements {\oe} with respect to the Moon are retrieved leveraging standard transformation schemes, e.g., in Vallado \cite{vallado2001fundamentals}, i.e., $\text{\oe} = \text{\oe}(\pvec{R}_E, \pvec{R}_E')$. The classical Keplerian elements are adopted in the current investigation as follows: $a$ (semi-major axis), $e$ (eccentricity), $i$ (inclination), $\omega$ (argument of perilune), $\Omega$ (right ascension of the ascending node), and $M$ (mean anomaly). To facilitate analysis for regions in the vicinity of the Moon, it is common practice to \textit{average} the perturbations from the Earth along the lunar orbit with respect to both $M$ and $\Omega$. This process results in the \textit{doubly-averaged} mean orbital elements, $\bar{\bar{\text{\oe}}}$, where the double bar indicates the quantity is averaged twice. 

\item \ul{Moon-Centered J2000 Inertial Frame} (\acrshort{mci}): the unit vectors $\hat{X} -\hat{Y} - \hat{Z}$ are defined consistent with the J2000 inertial frame. The nd spacecraft position vector in the frame is rendered as \cite{park2024assessment},
\begin{align}
    \label{eq:mci_brf} \vec{R} = \frac{l}{l_*}\bm{C} (\vec{\rho} - \vec{\rho}_M) = X\hat{X} + Y\hat{Y}+ Z\hat{Z}, \quad \bm{C} = \mb \hat{x} & \hat{y} & \hat{z} \me,
\end{align}
\end{enumerate}
where $l, l_*$ are the instantaneous and mean Earth-Moon distances, respectively. The unit vectors $\hat{x}-\hat{y}-\hat{z}$ are constructed instantaneously based on the numerical ephemerides DE440 \cite{park2021jpl}. 
Building upon these frames, the transformations between the various state representations are summarized in Fig. \ref{fig:frame_transformation}. The transformation between doubly-averaged mean elements, $\bar{\bar{\text{\oe}}}$, and osculating elements, \oe, is a non-trivial process in general, where a suitable transformation depends on the instantaneous configurations in the system. In the current investigation, however, $\bar{\bar{\text{\oe}}} = \text{\oe}$, i.e., both elements are treated equivalently. Such an approximation is justified from the fact that any semi-analytical transformation between $\bar{\bar{\text{\oe}}}$ and $\text{\oe}$ still requires further refinement. This work focuses on the refinement algorithm itself, alleviating the dependency on an accurate transformation between $\bar{\bar{\text{\oe}}}$ and $\text{\oe}$. The rest of the transformations in Fig. \ref{fig:frame_transformation} are provided more trivially via Eqs. \eqref{eq:eof_brf}-\eqref{eq:mci_brf}. Note that the velocities, i.e., $\pvec{R}_E', \pvec{\rho}', \pvec{R}'$, are interdependent via derivatives for these equations.

\subsection{\acrfull{dam}}
Within the Moon-dominated dynamical regime, 
several assumptions are introduced to supply the \acrshort{dam} that is employed in the current investigation. Firstly, the point-mass Earth is the only perturbing body. While it is possible to extend the \acrshort{dam} to include the perturbations from the nonspherical gravity from the Moon \cite{nie2018lunar}, such extensions are not considered here as the focus remains on the regions where the most dominant perturbations originate from the Earth \cite{ely2005stable}. Then, the Earth's perturbing potential is truncated leveraging the Legendre expansion as detailed by Longuski et al. \cite{longuski2022introduction}. Subsequently, the perturbation from the Earth is averaged over (1) $\Omega$ and (2) $M$ to supply the doubly-averaged, truncated Earth perturbing potential, 
\begin{align}
    \da{P} = \frac{1}{32}n_E\da{a}^2 \left[ (1+3\cos 2 \da{i} )(2 + 3\da{e}^2) + 30\da{e}^2 \sin^2 \da{i} \cos 2 \da{\omega} \right].
\end{align}
where $n_E$ is the dimensional mean angular rate for the Earth around the Moon. Recall that the doubly-averaged mean elements, i.e., $\da{\text{\oe}}$, are supplied within the \acrshort{eof}. From Lagrange's planetary equations, the mean orbital elements evolve at the following rates \cite{broucke2003long, nie2018lunar, folta2006lunar},
\begin{align}
    \label{eq:dadt} & \frac{d\da{a}}{dt} = 0 \\
    & \frac{d\da{e}}{dt} = \frac{15}{8}\frac{(1-\mu)n_E^2t_*}{\da{n}}\da{e}(1-\da{e}^2)^{1/2}\sin^2 \da{i}\sin 2 \da{\omega} \\
    & \frac{d\da{i}}{dt} = - \frac{15}{16}\frac{(1-\mu)n_E^2t_*}{\da{n}}\frac{\da{e}^2}{(1-\da{e}^2)^{1/2}}\sin 2\da{i} \sin 2\omega \\
    & \frac{d\da{\omega}}{dt} = \frac{3}{16}\frac{(1-\mu)n_E^2t_*}{\da{n}}\frac{1}{(1-\da{e}^2)^{1/2}} \left[ (3 + 2\da{e}^2 + 5\cos 2\da{i}) + 5(1-2\da{e}^2 - \cos 2\da{i} ) \cos 2\da{\omega} \right] \\
   \label{eq:dOmegadt}   & \frac{d\da{\Omega}}{dt} = \frac{3}{8}\frac{(1-\mu)n_E^2t_*}{\da{n}}\frac{1}{(1-\da{e}^2)^{1/2}} (5\da{e}^2 \cos 2\da{\omega} - 3\da{e}^2 - 2) \cos \da{i} \\
   \label{eq:dthetadt} &  \frac{d\da{M}}{dt} = \da{n}t_*,
\end{align}
where $\da{n}$ is the dimensional mean motion for the spacecraft (satellite) around the Moon supplied with $\da{a}$.

\subsection{\acrfull{cr3bp}}

The \acrshort{cr3bp} dynamics are cast within the BRF. The equations of motion are rendered as \cite{szebehely2012theory},
\begin{align}
    \label{eq:cr3bp} \pvec{\rho}'' = -2\hat{z}  \times {\pvec{\rho}}'+\nabla V.
\end{align}
The double prime indicates a second derivative with respect to the nd time, $t$. The vector cross product is denoted as $\times$, and $V$ is the pseudo-potential function rendered as $V = 1/2\cdot(x^2 + y^2) + (1-\mu)/d + \mu/r $ where $d$ and $r$ correspond to the nd distance between the spacecraft and Earth, and the spacecraft and Moon, respectively.

\subsection{\acrfull{hfem}}

The higher-fidelity modeling in the current investigation includes point-mass gravity from three celestial bodies: the Sun ($S$), Earth ($E$), and Moon ($M$). The instantneuos locations of the bodies are retrieved from the JPL ephemerides DE440 \cite{park2021jpl}.  With the Moon as the central gravitational body, the equations of motion for the spacecraft within the \acrshort{mci} are rendered,
\begin{align}
    \label{eq:hfem} \pvec{R}'' = -\frac{\mu}{R^3}\pvec{R} + \sum_{j = E, S} \left[ \frac{\mu_j}{R_{jM}
    ^3}\pvec{R}_{jM} - \frac{\mu_j}{R_{jc}
    ^3} \pvec{R}_{jc} \right],
\end{align}
where $\pvec{R}_{jc}$ and $\pvec{R}_{jM}$ denote the nd position vectors from the perturbing body $j$ to the spacecraft and the Moon, respectively. The nd gravitational parameter $\mu_j$ is defined as $\mu_j = \tilde{\mu}_j/(\tilde{\mu}_E + \tilde{\mu}_M) $. Subscripts $E, S$ are leveraged to indicate the Earth and the Sun, respectively. 

\subsection{\acrfull{qpos} and Frequency Decomposition}

\acrfull{qpos} arise naturally in the models incorporated throughout this investigation, but their characteristics depend on the underlying dynamical system. The \acrfull{dam} is an integrable Hamiltonian system with three degrees of freedom and three integrals of motion. All solutions exhibit (quasi-)periodic behavior. The \acrfull{cr3bp}, while non-integrable, retains one integral of motion in the \acrshort{brf} and admits quasi-periodic solutions within certain regimes. Both models are time-independent and supply QPOs with the maximum dimension of $\inter = 3$. Due to the low dimensional complexity, existing strategies such as GMOS are suitable for construction of \acrshort{qpos} for these simplified models. In contrast, the \acrfull{hfem} introduces additional time-dependent perturbations, generally increasing the number of fundamental frequencies from $\inter$ to $\ttl = \inter + \ex$. Under resonant interactions, $\ttl = \inter + \ex -1$. This work focuses on dynamical regimes where the QPOs persist and analog structures are well-defined within the HFEM. Each QPO is characterized by $\ttl$ fundamental frequencies. The full state vector on the QPO evolves as an infinite Fourier series comprised of linear combinations of these frequencies. Accurate recovery of these frequency components is central to the corrections strategy developed in this work.

%% file: text/quasi_periodic_orbits_freq.tex
\section{\label{sec:frequency_refinement}Frequency Refinement Methods}

While the existence of \acrfull{qpos} analytically guarantees that a trajectory on such a torus is represented as a Fourier series, numerically reconstructing this structure poses significant challenges. Since these trajectories are not available in closed form but are instead generated through sequential numerical integration, the \acrfull{dft} is employed in place of the \acrfull{cft}. However, the inherent limitations of the \acrshort{dft} hinder accurate frequency recovery. These challenges include, but are not limited to, spectral leakage, limited frequency resolution, and aliasing from components above the Nyquist threshold \cite{lyons1997understanding}.

To supply accurate frequency structures associated with quasi-periodic trajectories, the information from the \acrshort{dft} is \textit{refined}. To formalize the frequency refinement process, a scalar signal $\signalt(\tnd)$ is sampled along a quasi-periodic trajectory $\traj$. Typical choices for $\signalt$ include position or velocity components, distances to primary bodies, or angular quantities relative to a reference direction. The objective is to construct a quasi-periodic approximation $\signala(\tnd)$ of the form,
\begin{align} \label{eq:signal} \signalt(\tnd) \approx \signala (\tnd) = \ampA{0} + \sum_{\idxA=1}^{\nA} \ampA{\idxA} \cos(\freqA{\idxA} \tnd + \phaseA{\idxA}) = \ampA{0} + \sum_{\idxA=1}^{\nA} \signalaidxA (\tnd), \end{align} 
where $\ampA{\idxA}$, $\freqA{\idxA}$, and $\phaseA{\idxA}$ represent the amplitude, frequency, and phase of the $\idxA$-th quasi-periodic component. The number of terms $\nA$ is user-defined. Due to the discrete nature of the \acrshort{dft}, however, the identified peaks generally do not align exactly with the true frequency components $\freqA{\idxA}$. Consequently, the corresponding amplitude and phase estimates from the DFT are often insufficient, necessitating a refinement process. The quantities to be refined at each step are collected as,\begin{align} \label{eq:freeV} \freeVE_\idxA = \begin{bmatrix} \freqA{\idxA} \\ \ampA{\idxA} \\ \phaseA{\idxA} \end{bmatrix}. \end{align}
Two frequency refinement strategies are employed to estimate $\freeVE$: (i) \acrfull{laskar} method \cite{laskar1999introduction} and (ii) \acrfull{gomez} method \cite{gomez2010collocation}. Both methods define a set of algebraic constraints, and a differential corrections process is used to iteratively refine $\freeVE$ to satisfy these conditions. The notation introduced below largely follows the formulation in \cite{gomez2010collocation}. First, consider \acrshort{dft} on a quasi-periodic trajectory, $\traj$. The time-domain signal is constructed as,
\begin{align}
    & \signalt(\tidxT), \quad (0 \leq \idxT \leq \sampleSize-1) \\
    \label{eq:time_dft} & \tidxT =  \tndS\idxT \\
    & \tspan = \tndS \sampleSize,
\end{align}
where $\tndS$ is the sample spacing and $\tspan$ is the total time span. The sample size is denoted as $\sampleSize$. The Discrete Fourier Transform (DFT) of $\signalt$ is defined as,
\begin{align}
    \label{eq:dft} \dft{\signalt}{\genericFreq} \coloneq \frac{1}{\sampleSize} \sum_{\idxT=0}^{\sampleSize-1} \signalt(\tidxT) \han(\idxT) \exp{ -\imag  \genericFreq\tidxT},
\end{align}
where $\genericFreq$ is a general notation for a frequency. The frequencies in the current analysis are represented in the units of radians per respective unit time (not to be confused with cycles per unit time). The Hann window function of order two, i.e., $\han(\idxT) =2/3(1- \cos( 2\pi \idxT/\sampleSize) )^2$, is applied to alleviate spectral leakage \cite{lyons1997understanding}. While $\tspan$ and $\sampleSize$ serve as parameters for the \acrshort{dft} in Eq. \eqref{eq:dft}, they are assumed to be arbitrarily fixed numbers for notational brevity. In the \textit{standard} \acrshort{dft}, frequency bins are defined at $(\sampleSize/2+1)$ equally spaced points, corresponding to,
\begin{align}
    & \fndS = \frac{2\pi}{\tspan} \\
    \label{eq:bins} & \fidxF  =  \fndS \idxF = \frac{2\pi\idxF}{\tspan}, \quad (0\leq \idxF \leq N/2) ,
\end{align}
where $\sampleSize$ is assumed to be even (as is typical when using an \acrfull{fft}). The bin width $\fndS$ is solely a function of $\tspan$, where generally, investigating signals over a longer horizon time results in enhanced frequency resolutions. Evaluating the DFT at bin $\fidxF$ yields,
\begin{align}
    \label{eq:dft_signal_bin} \dft{\signalt}{\fidxF} & = \frac{1}{\sampleSize} \sum_{\idxT=0}^{\sampleSize-1} \signalt(\tidxT) \han(\idxT) \exp{ \frac{-2\pi \imag \idxT \idxF }{\sampleSize}} =  \frac{1}{\sampleSize} \sum_{\idxT=0}^{\sampleSize-1} \signalt(\tidxT) \han(\idxT) \exp{ -\imag \fidxF \tidxT} \\
    & = \frac{1}{\sampleSize}\sum_{\idxT=0}^{\sampleSize-1} \biggl(  \signalt(\tidxT) \han(\idxT)\bigl(\cos( \fidxF\tidxT) -\imag  \sin(\fidxF\tidxT) \bigr)  \biggr).
\end{align}
The imaginary number is denoted as $\imag$, i.e., $\imag^2 = -1$. Now, define the following quantities, 
\begin{align}
    \cosf (\fidxF) & = 2 \text{Re}(\signalf(\fidxF)) \\
    \sinf (\fidxF) & = -2 \text{Im}(\signalf(\fidxF)),
\end{align}
connected to the cosine (real) and sine (imag) components of $\signalf(\fidxF)$. From the standard \acrshort{dft} analysis on the signal from Eq. \eqref{eq:dft_signal_bin}, \textit{peaks} emerge that satisfy,
\begin{align}
    |\signalf(\fidxF)| \geq |\signalf(\fidxFm)|\quad \text{and}\quad |\signalf(\fidxF)| \geq |\signalf(\fidxFp)|.
\end{align}
For a quasi-periodic signal $\signalt$, prominent peaks typically emerge in the vicinity of the true frequency components. The $\idxA$-th peak from the bin frequency (Eq.\eqref{eq:dft_signal_bin}) is denoted as,
\begin{align}
    \freqP{\idxA},
\end{align}
and,
\begin{align}
    \label{eq:adjacent} \freqPPM{\idxA} = \freqP{\idxA} \pm \fndS,
\end{align}
corresponding to the adjacent bins to $\fidxF$. The initial guess for the free variables (Eq. \eqref{eq:freeV}) associated with the $\idxA$-th peak is supplied from the \acrshort{dft} as,
\begin{align}
    \freqA{\idxA} & \approx \freqP{\idxA} \\
    \ampA{\idxA} & \approx 2 |\dft{\signalt}{\freqP{\idxA}}| = \sqrt{(\cosf (\freqP{\idxA})) ^2 + (\sinf (\freqP{\idxA}))^2} \\
    \phaseA{\idxA} & \approx \arctan_2 (-\sinf (\freqP{\idxA}), \cosf (\freqP{\idxA})),
\end{align}
where $\arctan_2$ is the four-quadrant inverse tangent. Alternatively, an optimization problem may be formulated to locate $\signalaidxA = \ampA{\idxA} \cos(\freqA{\idxA} \tnd + \phaseA{\idxA}) $ from Eq. \eqref{eq:signal} that minimizes the error $e = \sum_{\idxT = 0}^{\sampleSize-1}(\signalt(\tidxT)- \signalaidxA(\tidxT))^2 $. The initial guess is subsequently refined using either the \acrfull{laskar} or \acrfull{gomez} methods that are described in the following sections. In both methods, the quasi-periodic representation $\signala$ (Eq. \eqref{eq:signal}) is sequentially refined for $\idxA \leq \nA $.

\subsection{\acrfull{laskar} Method \cite{laskar1999introduction}}

At the $\idxA$-th step, the \acrshort{laskar} method seeks to refine the frequency estimate $\freqA{\idxA}$ such that the magnitude of the DFT for the signal $\signalt$ is maximized at this frequency. This refinement condition is defined as,
\begin{align}
    \label{eq:laskar_condition}
    \freqA{\idxA} = \argmax{\freqAArg} \left| \signalf(\freqAArg) \right|.
\end{align}
The \textit{auxiliary} \acrshort{dft}, $\signalf(\freqAArg)$, is the continuous-frequency DFT of $\signalt$, constructed as,
\begin{align}
    \signalf(\freqAArg) 
    &= \frac{1}{\sampleSize} \sum_{\idxT=0}^{\sampleSize-1} \signalt(\tidxT)\, \han(\idxT)\, e^{-\imag \freqAArg \tidxT} \\
    &= \frac{1}{\sampleSize} \sum_{\idxT=0}^{\sampleSize-1} \signalt(\tidxT)\, \han(\idxT)\left[ \cos(\freqAArg \tidxT) - \imag \sin(\freqAArg \tidxT) \right].
\end{align}
Note that $\freqAArg$ represents a continuous frequency and may differ from the discretized frequency bins $\fidxF$ defined in the standard \acrshort{dft}. The magnitude of $\signalf(\freqAArg)$ is expressed as,
\begin{align}
    \left| \signalf(\freqAArg) \right| = \frac{1}{2} \sqrt{ \cosf(\freqAArg)^2 + \sinf(\freqAArg)^2 },
\end{align}
where the real and imaginary components are defined as,
\begin{align}
    \cosf(\freqAArg) &= \frac{2}{\sampleSize} \sum_{\idxT=0}^{\sampleSize-1} \signalt(\tidxT) \han(\idxT) \cos(\freqAArg \tidxT), \\
    \sinf(\freqAArg) &= \frac{2}{\sampleSize} \sum_{\idxT=0}^{\sampleSize-1} \signalt(\tidxT) \han(\idxT) \sin(\freqAArg \tidxT).
\end{align}
Then, at the maximum frequency from Eq. \eqref{eq:laskar_condition}, following first derivative results in zero,
\begin{align}
    \label{eq:laskar_f1} \textcolor{red}{\text{L-NAFF Constraint 1}} \quad
    \frac{d}{d \freqA{\idxA}} \left| \signalf(\freqA{\idxA}) \right| 
    = \frac{1}{2 \sqrt{ \cosf^2 + \sinf^2 }} \left( \cosf\, \frac{d \cosf}{d \freqA{\idxA}} + \sinf\, \frac{d \sinf}{d \freqA{\idxA}} \right) = 0.
\end{align}
The derivatives of the cosine and sine components are,
\begin{align}
    \frac{d \cosf}{d \freqA{\idxA}} &= -\frac{2}{\sampleSize} \sum_{\idxT=0}^{\sampleSize-1} \tidxT\, \signalt(\tidxT) \han(\idxT) \sin(\freqA{\idxA} \tidxT), \\
    \frac{d \sinf}{d \freqA{\idxA}} &= \frac{2}{\sampleSize} \sum_{\idxT=0}^{\sampleSize-1} \tidxT\, \signalt(\tidxT)\han(\idxT) \cos(\freqA{\idxA} \tidxT).
\end{align}
To refine the amplitude $\ampA{\idxA}$ and phase $\phaseA{\idxA}$, the $\idxA$-th approximated signal from Eq. \eqref{eq:signal}, i.e., 
\begin{align}
    \signalaidxA(\tnd) = \ampA{\idxA} \cos(\freqA{\idxA} \tnd + \phaseA{\idxA}),
\end{align}
is leveraged. The \acrshort{dft} on $\signalaidxA(\tnd)$ is required to match that of the original signal, i.e., 
\begin{align}
    \dft{\signalaidxA}{\freqA{\idxA}} = \dft{\signalt}{\freqA{\idxA}}.
\end{align}
Writing $\signalaidxA$ as $
    \signalaidxA(\tnd) = \ampA{\idxA} \cos\phaseA{\idxA} \cos(\freqA{\idxA} \tnd) - \ampA{\idxA} \sin\phaseA{\idxA} \sin(\freqA{\idxA} \tnd)$ and applying the DFT yields,
\begin{align}
    \dft{\signalaidxA}{\freqA{\idxA}} & = \frac{1}{2} \left[ \ampA{\idxA} \cos\phaseA{\idxA} \coscosf (\freqA{\idxA})- \ampA{\idxA} \sin\phaseA{\idxA} \cossinf (\freqA{\idxA}) \right] \nonumber \\ 
    &+ \frac{\imag}{2} \left[ -\ampA{\idxA} \cos\phaseA{\idxA} \sincosf(\freqA{\idxA}) + \ampA{\idxA} \sin\phaseA{\idxA} \sinsinf (\freqA{\idxA})\right],
\end{align}
where each auxiliary quantity is defined as:
\begin{align}
    \coscosf(\freqA{\idxA}) &= \frac{2}{\sampleSize} \sum_{\idxT=0}^{\sampleSize-1} \cos^2(\freqA{\idxA} \tidxT)\, \han(\idxT) \approx 1, \\
    \cossinf(\freqA{\idxA}) &= \frac{2}{\sampleSize} \sum_{\idxT=0}^{\sampleSize-1} \cos(\freqA{\idxA} \tidxT) \sin(\freqA{\idxA} \tidxT)\, \han(\idxT)  \approx 0, \\
    \sincosf(\freqA{\idxA}) &= \cossinf(\freqA{\idxA}) , \\
    \sinsinf(\freqA{\idxA}) &= \frac{2}{\sampleSize} \sum_{\idxT=0}^{\sampleSize-1} \sin^2(\freqA{\idxA} \tidxT)\, \han(\idxT)  \approx 1.
\end{align}
Matching the real and imaginary parts yields two additional constraints,
\begin{align}
    \label{eq:laskar_f2}&\textcolor{red}{\text{L-NAFF Constraint 2}} \quad
    \ampA{\idxA} \cos\phaseA{\idxA} \coscosf(\freqA{\idxA}) - \ampA{\idxA} \sin\phaseA{\idxA} \cossinf(\freqA{\idxA}) = \cosf, \\
    \label{eq:laskar_f3}&\textcolor{red}{\text{L-NAFF Constraint 3}}\quad
    \ampA{\idxA} \cos\phaseA{\idxA} \sincosf(\freqA{\idxA}) - \ampA{\idxA} \sin\phaseA{\idxA} \sinsinf(\freqA{\idxA}) = \sinf.
\end{align}
Together, Eqs. \eqref{eq:laskar_f1}, \eqref{eq:laskar_f2}, \eqref{eq:laskar_f3} define the constraint vector for \acrshort{laskar} as,
\begin{align}
    \label{eq:laskar_constraint} \constraintVL = \begin{bmatrix}
        \constraintC{L, \freqAArg} \\ \constraintC{L, \mathcal{C}} \\\constraintC{L, \mathcal{S}}  
    \end{bmatrix} = \begin{bmatrix}
        \frac{1}{2\sqrt{\cosf(\freqA{\idxA})^2 + \sinf(\freqA{\idxA})^2}}\left(\cosf(\freqA{\idxA}) \frac{d \cosf(\freqA{\idxA})}{d \freqA{\idxA}} + \sinf(\freqA{\idxA}) \frac{d \sinf(\freqA{\idxA})}{d \freqA{\idxA}}  \right) \\
        \ampA{\idxA}\cos \phaseA{\idxA}\coscosf(\freqA{\idxA})  - \ampA{\idxA}\sin \phaseA{\idxA} \cossinf(\freqA{\idxA})  - \cosf(\freqA{\idxA}) \\
    \ampA{\idxA}\cos \phaseA{\idxA}\sincosf(\freqA{\idxA})   - \ampA{\idxA}\sin \phaseA{\idxA} \sinsinf(\freqA{\idxA}) - \sinf(\freqA{\idxA}).
    \end{bmatrix} = \vec{0},
\end{align}
solved via differential corrections. The Jacobian matrix of $\constraintVL$ with respect to $\freeVE_\idxA$ renders,
\begin{align}
    \label{eq:laskar_jacobian}    \frac{\partial \constraintVL}{\partial \freeVE_\idxA} = \begin{bmatrix}
    \frac{\partial \constraintC{L, \freqAArg}}{\partial \freqA{\idxA}} & \frac{\partial \constraintC{L, \freqAArg}}{\partial \ampA{\idxA}} & \frac{\partial \constraintC{L, \freqAArg}}{\partial \phaseA{\idxA}} \\
        \frac{\partial \constraintC{L, \mathcal{C}}}{\partial \freqA{\idxA}} & \frac{\partial \constraintC{L, \mathcal{C}}}{\partial \ampA{\idxA}} & \frac{\partial \constraintC{L, \mathcal{C}}}{\partial \phaseA{\idxA}} \\
        \frac{\partial \constraintC{L, \mathcal{S}}}{\partial \freqA{\idxA}} & \frac{\partial \constraintC{L, \mathcal{S}}}{\partial \ampA{\idxA}} & \frac{\partial \constraintC{L, \mathcal{S}}}{\partial \phaseA{\idxA}}
    \end{bmatrix},
\end{align}
with each component derived analytically using the chain rule and second derivatives of $\cosf$ and $\sinf$, included in \hyperref[ap:jacobian]{Appendix A}. After convergence of the $\idxA$-th frequency component, the next peak $(\idxA + 1)$ is detected from the residual signal,
\begin{align}
    \signalt^{new}_{\idxA} \coloneq \signalt - (\ampA{0} + \sum_{a=1}^{\idxA }\tilde{q}_a)
\end{align}
In the new signal $\signalt^{new}_{\idxA}$, the $(\idxA+1)$-th peak from the original signal $\signalt$ is expected to occur as the most dominant peak, facilitating the isolation of subsequent lower-amplitude components. The zeroth-order term $\ampA{0}$ is trivially constructed as the average of $\signalt$.

\subsection{\acrfull{gomez} Method \cite{gomez2010collocation}}

In contrast to the \acrshort{laskar} method that directly maximizes the DFT magnitude at the approximated frequency, the \acrshort{gomez} approach refines the frequency estimate $\freqA{\idxA}$ by enforcing collocation conditions in adjacent frequency bins. This strategy relies on matching the DFT of the approximated signal $\signalaidxA$ to that of the original signal $\signalt$ at these selected frequencies. To formulate the algebraic constraints, consider the DFTs of cosine and sine functions evaluated at $\fidxF$, the $k$-th bin from the \textit{standard} \acrshort{dft},
\begin{align}
    \signalcosf(\fidxF) &= \frac{1}{\sampleSize} \sum_{\idxT=0}^{\sampleSize-1} \cos(\freqA{\idxA} \tidxT)\, \han(\idxT)\, e^{-\imag \fidxF \tidxT}, \\
    \signalsinf(\fidxF) &= \frac{1}{\sampleSize} \sum_{\idxT=0}^{\sampleSize-1} \sin(\freqA{\idxA} \tidxT)\, \han(\idxT)\, e^{-\imag \fidxF \tidxT}, \\
    \frac{d \signalcosf(\fidxF)}{d \freqA{\idxA}} &= -\frac{1}{\sampleSize} \sum_{\idxT=0}^{\sampleSize-1} \tidxT \sin(\freqA{\idxA} \tidxT)\, \han(\idxT)\, e^{-\imag \fidxF \tidxT}, \\
    \frac{d \signalsinf(\fidxF)}{d \freqA{\idxA}} &= \frac{1}{\sampleSize} \sum_{\idxT=0}^{\sampleSize-1} \tidxT \cos(\freqA{\idxA} \tidxT)\, \han(\idxT)\, e^{-\imag \fidxF \tidxT}.
\end{align}
The following real-valued quantities are defined for use in the constraint expressions,
\begin{align}
    \coscosf(\fidxF) &= 2\, \text{Re}(\signalcosf(\fidxF)), & \quad
    \cossinf(\fidxF) &= 2\, \text{Re}(\signalsinf(\fidxF)), \\
    \sincosf(\fidxF) &= -2\, \text{Im}(\signalcosf(\fidxF)), & \quad
    \sinsinf(\fidxF) &= -2\, \text{Im}(\signalsinf(\fidxF)), \\
    \frac{d \coscosf(\fidxF)}{d \freqA{\idxA}} &= 2\, \text{Re} \left( \frac{d \signalcosf(\fidxF)}{d \freqA{\idxA}} \right), & \quad
    \frac{d \cossinf(\fidxF)}{d \freqA{\idxA}} &= 2\, \text{Re} \left( \frac{d \signalsinf(\fidxF)}{d \freqA{\idxA}} \right), \\
    \frac{d \sincosf(\fidxF)}{d \freqA{\idxA}} &= -2\, \text{Im} \left( \frac{d \signalcosf(\fidxF)}{d \freqA{\idxA}} \right), & \quad
    \frac{d \sinsinf(\fidxF)}{d \freqA{\idxA}} &= -2\, \text{Im} \left( \frac{d \signalsinf(\fidxF)}{d \freqA{\idxA}} \right).
\end{align}
Based on this auxiliary \acrshort{dft} information, the \acrshort{gomez} refinement process enforces collocation at two frequencies: (1) the $\idxA$-th DFT peak, $\freqP{\idxA}$, and (2) one of its adjacent bins, $\freqPPM{\idxA}$ from Eq. \eqref{eq:adjacent}. Equivalently, 
\begin{align}
    \label{eq:gomez_match1}
    \dft{\signalaidxA}{\freqP{\idxA}} &= \dft{\signalt}{\freqP{\idxA}}, \\
    \label{eq:gomez_match2}
    \dft{\signalaidxA}{\freqPO{\idxA}} &= \dft{\signalt}{\freqPO{\idxA}}.
\end{align}
The auxiliary frequency $\freqPO{\idxA}$ is selected from $\freqPM{\idxA}$ and $\freqPP{\idxA}$ based on the larger DFT amplitude of the original signal. Specifically, $\freqPO{\idxA}$ corresponds to the neighboring bin that yields the greater value of $|\dft{\signalt}{\freqPO{\idxA}}|$. Since each DFT match from Eqs. \eqref{eq:gomez_match1}-\eqref{eq:gomez_match2} introduces two scalar conditions (real and imaginary parts), these constraints overdetermine the three-element state vector $\freeVE_\idxA$ from Eq. \eqref{eq:freeV}. To ensure an invertible Jacobian, only one component (cosine or sine) is enforced at $\freqPO{\idxA}$, selected based on conditioning \cite{gomez2010collocation2}. The resulting three scalar constraints are:
\begin{align}
    \label{eq:gomez_constraint1} & \textcolor{red}{\text{GMS-C Constraint 1}} \quad 
    \ampA{\idxA} \cos\phaseA{\idxA}\, \coscosf(\freqP{\idxA}) - \ampA{\idxA} \sin\phaseA{\idxA}\, \cossinf(\freqP{\idxA}) = \cosf(\freqP{\idxA}), \\
    \label{eq:gomez_constraint2} & \textcolor{red}{\text{GMS-C Constraint 2}}\quad 
    \ampA{\idxA} \cos\phaseA{\idxA}\, \sincosf(\freqP{\idxA}) - \ampA{\idxA} \sin\phaseA{\idxA}\, \sinsinf(\freqP{\idxA}) = \sinf(\freqP{\idxA}), \\
    \label{eq:gomez_constraint3} & \textcolor{red}{\text{GMS-C Constraint 3}}\quad 
    \ampA{\idxA} \cos\phaseA{\idxA}\, \cscosf(\freqPO{\idxA}) - \ampA{\idxA} \sin\phaseA{\idxA}\, \cssinf(\freqPO{\idxA}) = \csf(\freqPO{\idxA}),
\end{align}
where $\mathcal{CS}$ denotes that one of the $\mathcal{C}$ or $\mathcal{S}$ is selected based on the condition from \hyperref[ap:jacobian]{Appendix A}. Three collocation conditions (Eqs. \eqref{eq:gomez_constraint1}-\eqref{eq:gomez_constraint3}) are concatenated as,
\begin{align}
    \label{eq:gomez_constraint_idxA} \constraintVG = \begin{bmatrix}
        \constraintC{G, \mathcal{C}(\freqP{\idxA})} \\ 
        \constraintC{G, \mathcal{S}(\freqP{\idxA})} \\
        \constraintC{G, \mathcal{CS}(\freqPO{\idxA})}
    \end{bmatrix} = \begin{bmatrix}
        \ampA{\idxA}\cos\phaseA{\idxA} \coscosf(\freqP{\idxA}) - \ampA{\idxA}\sin\phaseA{\idxA} \cossinf(\freqP{\idxA}) - \cosf(\freqP{\idxA}) \\
        \ampA{\idxA}\cos\phaseA{\idxA} \sincosf(\freqP{\idxA}) - \ampA{\idxA}\sin\phaseA{\idxA} \sinsinf(\freqP{\idxA}) - \sinf(\freqP{\idxA}) \\
        \ampA{\idxA}\cos\phaseA{\idxA} \cscosf(\freqPO{\idxA}) - \ampA{\idxA}\sin\phaseA{\idxA} \cssinf(\freqPO{\idxA}) - \csf(\freqPO{\idxA})
    \end{bmatrix} = \vec{0}.
\end{align}
The Jacobian $\partial \constraintVG / \partial \freeVE$ is constructed as follows,
\begin{align}
    \label{eq:gomez_jacobian}    \frac{\partial \constraintVG}{\partial \freeVE_\idxA} = \begin{bmatrix}
    \frac{\partial \constraintC{G, \mathcal{C}(\freqP{\idxA})}}{\partial \freqA{\idxA}} & \frac{\partial \constraintC{G, \mathcal{C}(\freqP{\idxA})}}{\partial \ampA{\idxA}} & \frac{\partial \constraintC{G, \mathcal{C}(\freqP{\idxA})}}{\partial \phaseA{\idxA}} \\
        \frac{\partial \constraintC{G, \mathcal{S}(\freqP{\idxA})}}{\partial \freqA{\idxA}} & \frac{\partial \constraintC{G, \mathcal{S}(\freqP{\idxA})}}{\partial \ampA{\idxA}} & \frac{\partial \constraintC{G, \mathcal{S}(\freqP{\idxA})}}{\partial \phaseA{\idxA}} \\
        \frac{\partial \constraintC{G, \mathcal{CS}(\freqPO{\idxA})}}{\partial \freqA{\idxA}} & \frac{\partial \constraintC{G, \mathcal{CS}(\freqPO{\idxA})}}{\partial \ampA{\idxA}} & \frac{\partial \constraintC{G, \mathcal{CS}(\freqPO{\idxA})}}{\partial \phaseA{\idxA}}
    \end{bmatrix},
\end{align}
where each component is supplied in \hyperref[ap:jacobian]{Appendix A}.

\subsection{Numerical Comparison of the Two Refinement Strategies}

The two refinement strategies, \acrshort{laskar} and \acrshort{gomez}, share a 
common goal to refine quasi-periodic approximations with auxiliary \acrshort{dft}-based constraints. The auxiliary DFTs yield algebraic conditions (Eqs. \eqref{eq:laskar_constraint} and \eqref{eq:gomez_constraint_idxA}) 
that are enforced by adjusting the free-variable vector 
(Eq.~\eqref{eq:freeV}). Recall that the \textit{standard} \acrshort{dft} refers to the signal transform 
in Eq. \eqref{eq:dft} evaluated at the discrete frequency bins 
(Eq. \eqref{eq:bins}). The two refinement methods differ then as follows: (1) \acrshort{laskar} applies the DFT directly to $\signalt$, but at continuous frequencies $\freqA{\idxA}$ that generally do not align with the bin centers (i.e., $\freqA{\idxA} \neq \fidxF$), and (2) \acrshort{gomez} instead evaluates the DFT at the standard bins, but on a cosine and sine basis signals defined by $\freqA{\idxA}$. Another distinction emerges when examining the Jacobian matrices in Eqs.~\eqref{eq:laskar_jacobian} and \eqref{eq:gomez_jacobian}. In the \acrshort{laskar} strategy, the constraint vector at the $\idxA$-th step depends only on the corresponding free-variable vector $\freeVE_\idxA$. In contrast, the \acrshort{gomez} method introduces interdependencies across modes: the Jacobian $\partial \constraintV{l} / \partial \freeVE_\idxA$ is generally nonzero for all $l \neq \idxA$. As a result, the \acrshort{gomez} strategy refines the first $\idxA$ frequency components simultaneously at the $\idxA$-th step. Moreover, the zeroth-order coefficient $\ampA{0}$ is also subject to refinement in the \acrshort{gomez} formulation. While this aspect is not directly leveraged in the current investigation, further details are found in G{\'o}mez, Mondelo, and Sim{\'o}~\cite{gomez2010collocation}.

\begin{figure}[htpb]
    \centering
    \includegraphics[width=0.5\linewidth]{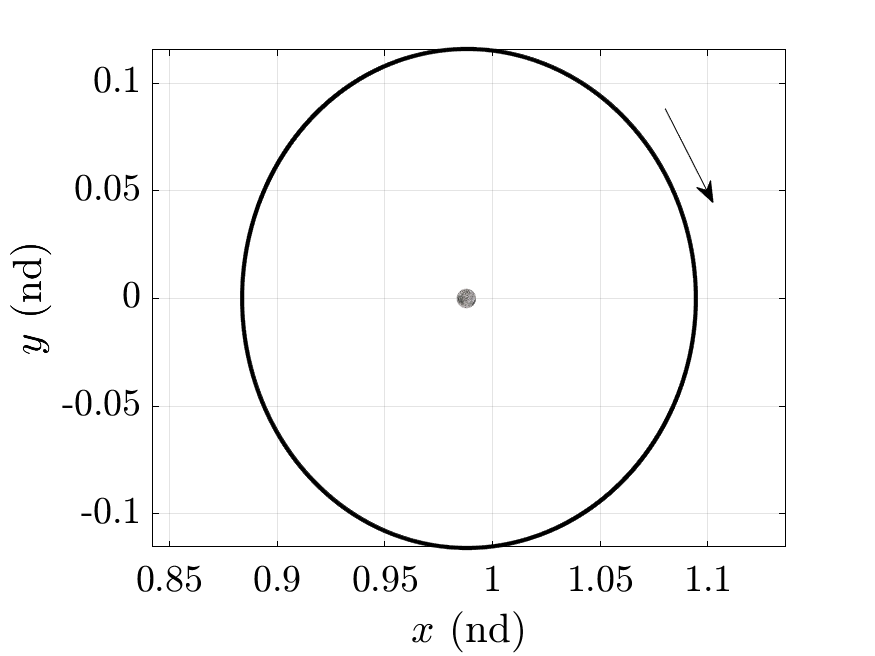}
    \caption{Sample \acrshort{cr3bp} periodic-\acrshort{dro} within the \acrfull{brf}}
    \label{fig:dro-compare-l-g-orbit}
\end{figure}

\begin{figure}[htpb]
    \centering
    \includegraphics[width=0.8\linewidth]{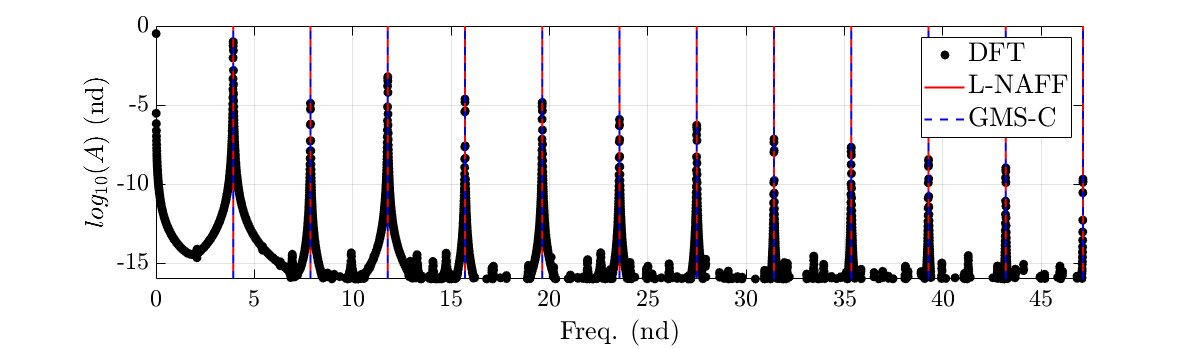}
    \caption{\acrshort{dft} on $x$-position components for the periodic-\acrshort{dro} in Fig. \ref{fig:dro-compare-l-g-orbit}}
    \label{fig:dro-compare-l-g-dft.png}
\end{figure}

\begin{figure}[htpb]
    \centering
    \includegraphics[width=0.8\linewidth]{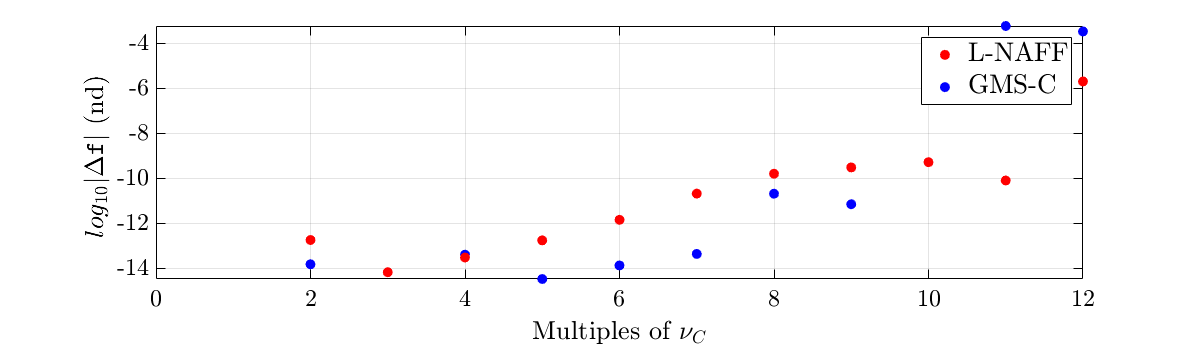}
    \caption{Errors for refined frequencies at the multiples of $\fund{C}$}
    \label{fig:dro-compare-l-g-refine.png}
\end{figure}

To illustrate the difference in performance between these two strategies, consider a sample \acrshort{cr3bp} \acrfull{dro} as illustrated in Fig. \ref{fig:dro-compare-l-g-orbit}. The initial condition for the orbit is supplied in Table \ref{table:dro_initial}. Along the periodic orbit, $\sampleSize = 2^{16}$ samples on the $x$-position components are constructed and leveraged as $\signalt$ for the \acrshort{dft} over 20 years. Note that the period of the \acrshort{dro} is $1.6$ (nd), where the fundamental frequency corresponds to $\fund{C} = 2\pi/1.6 \approx 3.927$ (nd). The resulting DFT spectrum is included in Fig.~\ref{fig:dro-compare-l-g-dft.png}, with dot markers at each frequency bin. The vertical axis illustrates $\log_{10}$ of the DFT amplitude ($\ampA{}$). While the peaks are generally associated with multiples of $\fund{C}$, due to the aforementioned limitations in the \acrshort{dft}, refinement strategies are employed. The results are illustrated in red and blue within Fig. \ref{fig:dro-compare-l-g-dft.png} for the \acrshort{laskar} and \acrfull{gomez} methods, respectively. Figure~\ref{fig:dro-compare-l-g-refine.png} illustrates the frequency 
errors (in comparison to integer multiples of $\fund{C}$) 
after applying each refinement. No marker appears for
$\fund{C}$ (multiple 1), as both methods resolve 
it to near-zero error. While the detailed performance varies depending on factors like numerical formulation and sampling resolution, \textbf{both methods supply accurate estimates for the dominant peaks within the lower-frequeny domain}. Throughout the current investigation, both strategies are leveraged interchangeably.

\begin{table}[htpb]
\centering
\caption{\label{table:dro_initial}Period and initial state for the periodic-\acrshort{dro} represented in Fig. \ref{fig:dro-compare-l-g-orbit}}
\begin{tabular}{|c|c|}
\hline
Variable       & Value (nd) \\ \hline
period &  $1.6$     \\ \hline
$x$       &    $0.883749964899239$    \\ \hline
$y$    &    $0$   \\ \hline
$z$    &   $0$    \\ \hline
$x'$    &  $0$     \\ \hline
$y'$   &    $0.470425740470053  $   \\ \hline
      $z'$&   $0$    \\ \hline
\end{tabular}
\end{table}

%% file: text/single_shooter.tex
\section{\label{sec:single}\acrfull{fdc}: Single-Shooting Formulation}

The overarching goal of this investigation is to develop a \acrfull{fdc} capable of generating $\ttl$-dimensional \acrfull{qpos}. The preceding section introduces two refinement strategies that, given a fixed signal $\signalt$, extract accurate representations of individual frequency components, i.e., $\freqA{\idxA}$, $\ampA{\idxA}$, and $\phaseA{\idxA}$, associated with the $\idxA$-th peak in the frequency-domain. In contrast, the \acrshort{fdc} adjusts the signal itself so that these spectral components match desired target values, thereby controlling the underlying $\ttl$-dimensional QPO. 

This section (single-shooting formulation) focuses on dynamical regimes that are locally stable, i.e., those free of strong instabilities. In such settings, the signal $\signalt$ is supplied by propagating a single six-dimensional initial state $\x{0} \in \real{6}$ over a long time-span. To establish the core algorithm, a generic, autonomous $\inter$ degree-of-freedom Hamiltonian system is examined as the governing dynamics. Within such systems, dynamical regimes may exist where dense $\inter$-dimensional QPOs exist locally, demonstrating a stable behavior.\footnote{Although resonant periodic orbits in the vicinity of the dense $\inter$-dimensional QPOs with $\inter > 2$ may introduce weak instabilities due to Arnold diffusion \cite{arnol2013mathematical}, their growth rates are super-exponentially slow \cite{morbidelli1995superexponential}.} Relevant examples in the Earth–Moon \acrshort{cr3bp} include, but are not limited to, the \acrfull{dro} family and the \acrfull{elfos}. The single-shooting \acrshort{fdc} formulation assumes the initial state $\x{0}$ as the free variable and enforces the following constraint, 
\begin{align}
    \label{eq:frequency_constraint} \f{f} = \begin{bmatrix} 
        \freqA{\idxA} - \freqA{\idxA}^T \\
        \phaseA{\idxA} - \phaseA{\idxA}^T
    \end{bmatrix}
    \quad \text{or} \quad
    \begin{bmatrix}
        \ampA{\idxA} - \ampA{\idxA}^T \\
        \phaseA{\idxA} - \phaseA{\idxA}^T
    \end{bmatrix},
\end{align}
where the superscript $T$ denotes the target frequency-domain quantities. Importantly, the amplitude and frequency components $\ampA{\idxA}$ and $\freqA{\idxA}$ are generally coupled within a QPO. Attempting to target both simultaneously typically results in an over-constrained system. Therefore, the \acrshort{fdc} formulation targets either $\freqA{\idxA}$ or $\ampA{\idxA}$, depending on the application. To construct a differential corrector, the gradient of $\f{f}$ with respect to $\x{0}$ is required. Specifically, the following 
sensitivities must be supplied,
\begin{align}
    \label{eq:fdt_partial_gomez_required}
    \frac{\partial \freqA{\idxA}}{\partial \x{0}}, \quad 
    \frac{\partial \ampA{\idxA}}{\partial \x{0}}, \quad 
    \frac{\partial \phaseA{\idxA}}{\partial \x{0}}.
\end{align}
These derivatives are constructed by differentiating the algebraic constraints associated with the refinement strategies introduced from Eqs.~\eqref{eq:laskar_constraint} and \eqref{eq:gomez_constraint_idxA}. The required sensitivity expressions for each method are presented in the following subsections.

\subsection{\label{sec:sensitivity}Sensitivities of the State Vector on the Frequency Components }

The sensitivities of $\freqA{\idxA}$, $\ampA{\idxA}$, and $\phaseA{\idxA}$ with respect to the initial state $\x{0}$ are derived for both refinement methods. This section focuses on the single-shooting formulation, where the signal $\signalt$ is generated by forward propagation of $\x{0}$ under any dynamical model. Note that,
\begin{align}
    \label{eq:xp}\xp{0}{\tidxT}{\x{0}}{\epoch{0}}
\end{align}
denotes the state propagated from $\x{0}$ initiating at the epoch $\epoch{0}$ (assumed to be fixed), propagated for $\tidxT$. Then, the sensitivity of the propagated state as the initial state evolves is captured by the State Transition Matrix (STM). The STM from the initial state $\x{0}$ at epoch $t_0$ to time $\tidxT$ is denoted $\stm{0}{\tidxT}{\x{0}}{\epoch{0}}$, and satisfies,
\begin{align}
    \frac{d \xp{0}{\tidxT}{\x{0}}{\epoch{0}}}{d \x{0}} = \stm{0}{\tidxT}{\x{0}}{\epoch{0}}.
\end{align}
The sensitivity of any scalar signal $\signalt(\tidxT)$ derived from the state vector is produced through the chain rule,
\begin{align}
    \label{eq:dsignal_dx}\frac{d \signalt(\tidxT)}{d \x{0}} = 
    \frac{d \signalt(\tidxT)}{d \xp{0}{\tidxT}{\x{0}}{\epoch{0}}} 
    \cdot\stm{0}{\tidxT}{\x{0}}{\epoch{0}}.
\end{align}
Depending on the model and the frame, $\signalt$ may or may not correspond trivially to a state component from $\xpp$, illustrated further in the subsequent numerical examples. 

\subsubsection{Sensitivities within the \acrfull{laskar} Method}

The sensitivities in Eq. \eqref{eq:fdt_partial_gomez_required} are characterized via differentiating the constraint vector from the \acrshort{laskar} method, i.e., Eq. \eqref{eq:laskar_constraint}. The process yields, 
\begin{align}
     \label{eq:fdt_partial_laskar} \frac{\partial \constraintVL }{\partial \freeVE_\idxA} \mb \frac{\partial\freqA{\idxA}}{\partial\x{0}} \\ \frac{\partial\ampA{\idxA}}{\partial\x{0}} \\ \frac{\partial\phaseA{\idxA}}{\partial\x{0}} \me = -\mb \frac{\partial \constraintC{L, \freqA{}} }{\partial \x{0}}\\ \frac{\partial \constraintC{L, \mathcal{C}} }{\partial \x{0}}\\ \frac{\partial \constraintC{L, \mathcal{S}} }{\partial \x{0}}\me.
\end{align}
The first term, $\frac{\partial \constraintVL }{\partial \freeVE_\idxA} $, is an invertible matrix constructed from Eq. \eqref{eq:laskar_jacobian}. The vector components on the right side of Eq. \eqref{eq:fdt_partial_laskar} involve expressions from Eq. \eqref{eq:dsignal_dx} and are supplied in \hyperref[ap:sensitivities]{Appendix B}. Then, it is possible to invert the matrix and write,
\begin{align}
      \label{eq:sensitivity_laskar} \mb \frac{\partial\freqA{\idxA}}{\partial\x{0}} \\ \frac{\partial\ampA{\idxA}}{\partial\x{0}} \\ \frac{\partial\phaseA{\idxA}}{\partial\x{0}} \me = - \left(\frac{\partial \constraintVL }{\partial \freeVE_\idxA }\right)^{-1} \mb \frac{\partial \constraintC{L, \freqA{}} }{\partial \x{0}}\\ \frac{\partial \constraintC{L, \mathcal{C}} }{\partial \x{0}}\\ \frac{\partial \constraintC{L, \mathcal{S}} }{\partial \x{0}}\me, 
\end{align}
supplying the sensitivities in Eq. \eqref{eq:fdt_partial_gomez_required}. Note that the left side of Eq. \eqref{eq:sensitivity_laskar} is a $(3\times 6)$ matrix.

\subsubsection{Sensitivities within the \acrfull{gomez} Method}

A similar process is now initiated with the \acrshort{gomez} method. Differentiating Eq. \eqref{eq:gomez_constraint_idxA} yields,
\begin{align}
     \label{eq:fdt_partial_gomez} \frac{\partial \constraintVG }{\partial \freeVE} \mb \frac{\partial\freqA{\idxA}}{\partial\x{0}} \\ \frac{\partial\ampA{\idxA}}{\partial\x{0}} \\ \frac{\partial\phaseA{\idxA}}{\partial\x{0}} \me = -\mb \frac{\partial  \constraintC{G, \mathcal{C}(\freqP{\idxA})}}{\partial \x{0}} \\ \frac{\partial  \constraintC{G, \mathcal{S}(\freqP{\idxA})}}{\partial \x{0}} \\ \frac{\partial  \constraintC{G, \mathcal{CS}(\freqP{\idxA})}}{\partial \x{0}} \me = \mb \frac{\partial\cosf(\freqP{\idxA})}{\partial\x{0}} \\ \frac{\partial\sinf(\freqP{\idxA})}{\partial\x{0}} \\ \frac{d\csf(\freqPO{\idxA})}{\partial\x{0}} \me, 
\end{align}
where $\frac{\partial \constraintVG }{\partial \freeVE_\idxA}$ is an invertible $3 \times 3$ matrix supplied via Eq. \eqref{eq:gomez_jacobian}. The expressions on the right side of Eq. \eqref{eq:fdt_partial_gomez} involve expressions from Eq. \eqref{eq:dsignal_dx} and are supplied in \hyperref[ap:sensitivities]{Appendix B}. Then, it is possible to invert the matrix and write,
\begin{align}
    \mb \frac{\partial\freqA{\idxA}}{\partial\x{0}} \\ \frac{\partial\ampA{\idxA}}{\partial\x{0}} \\ \frac{\partial\phaseA{\idxA}}{\partial\x{0}} \me = \left(  \frac{\partial \constraintVG }{\partial \freeVE}\right)^{-1}\mb \frac{\partial\cosf(\freqP{\idxA})}{\partial\x{0}} \\ \frac{\partial\sinf(\freqP{\idxA})}{\partial\x{0}} \\ \frac{d\csf(\freqPO{\idxA})}{\partial\x{0}} \me,
\end{align}
deriving the sensitivities for the \acrshort{gomez} approach.

\subsection{Applications}

Thus far, the \acrfull{fdc} is derived within a general single-shooting formulation. The general aspects of the algorithm are further illustrated in the following examples focusing on the \acrfull{dros} and \acrfull{elfos}.

\subsubsection{Specifying a Quasi-Periodic Trajectory within the \acrfull{dro} Family}

\begin{figure}[h!]
    \centering
    \begin{subfigure}[b]{0.48\textwidth}
        \centering
        \includegraphics[width = 0.99\textwidth]{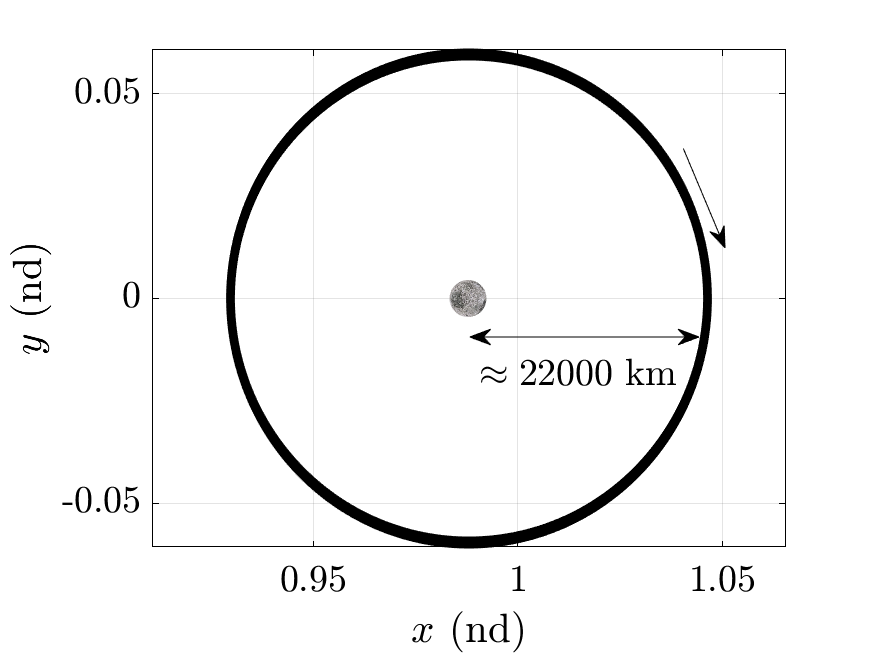}
        \caption{Geometry within the \acrshort{brf}}
        \label{fig:dro-cr3bp-geometry}
    \end{subfigure}
    \begin{subfigure}[b]{0.48\textwidth}
        \centering
        \includegraphics[width = 0.99\textwidth]{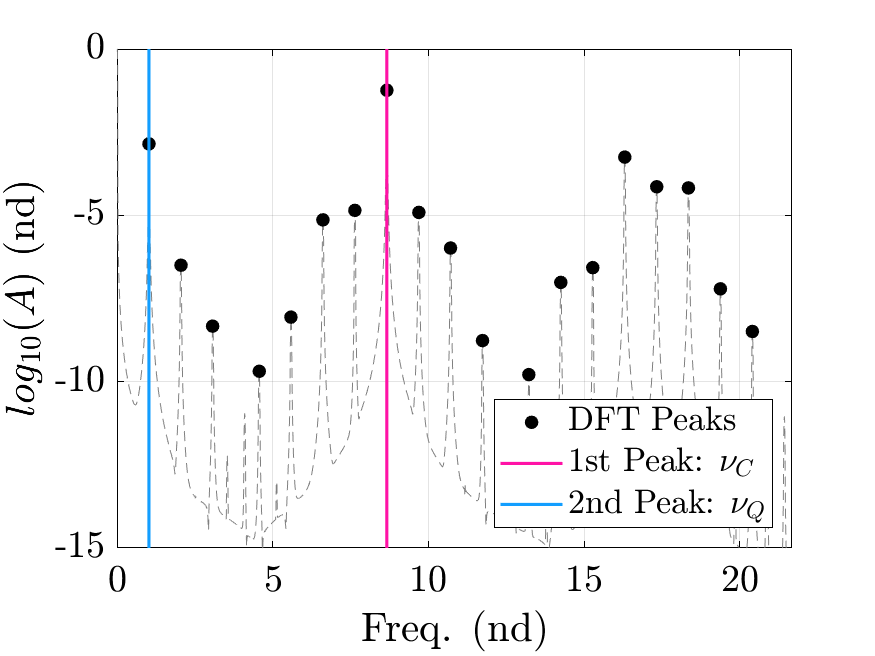}
        \caption{\acrshort{dft} and refined peaks with $\signalt = x$ (\acrshort{laskar})}
        \label{fig:dro-cr3bp-dft}
    \end{subfigure}
  \caption{\acrshort{cr3bp} quasi-\acrshort{dro} initiated with the state in Table \ref{table:dro-ig}}
  \label{fig:dro-cr3bp}
\end{figure}

\begin{figure}[h!]
    \centering
    \begin{subfigure}[b]{0.48\textwidth}
        \centering
        \includegraphics[width = 0.99\textwidth]{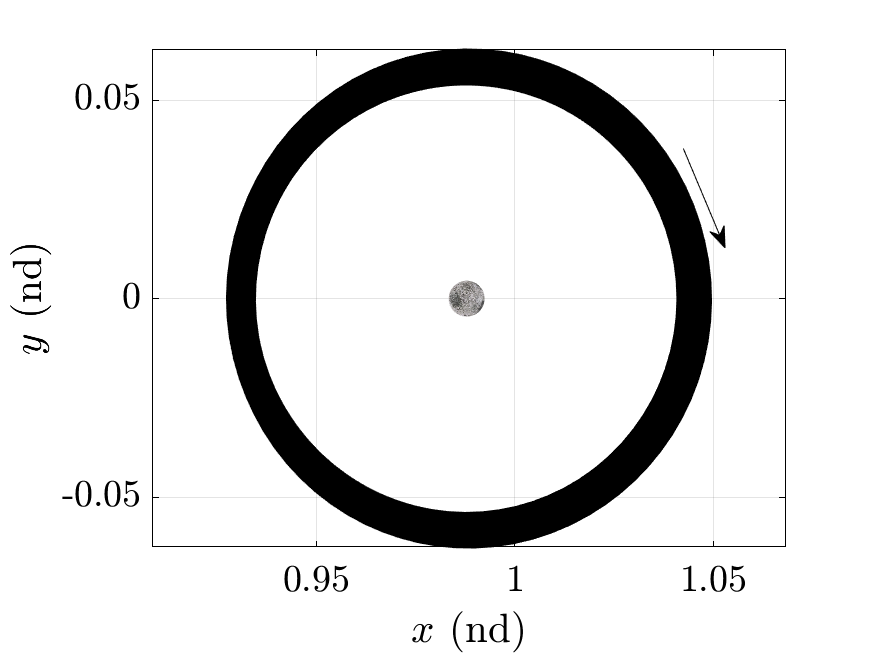}
        \caption{Geometry within the \acrshort{brf}}
        \label{fig:dro-hfem-geometry}
    \end{subfigure}
    \begin{subfigure}[b]{0.48\textwidth}
        \centering
        \includegraphics[width = 0.99\textwidth]{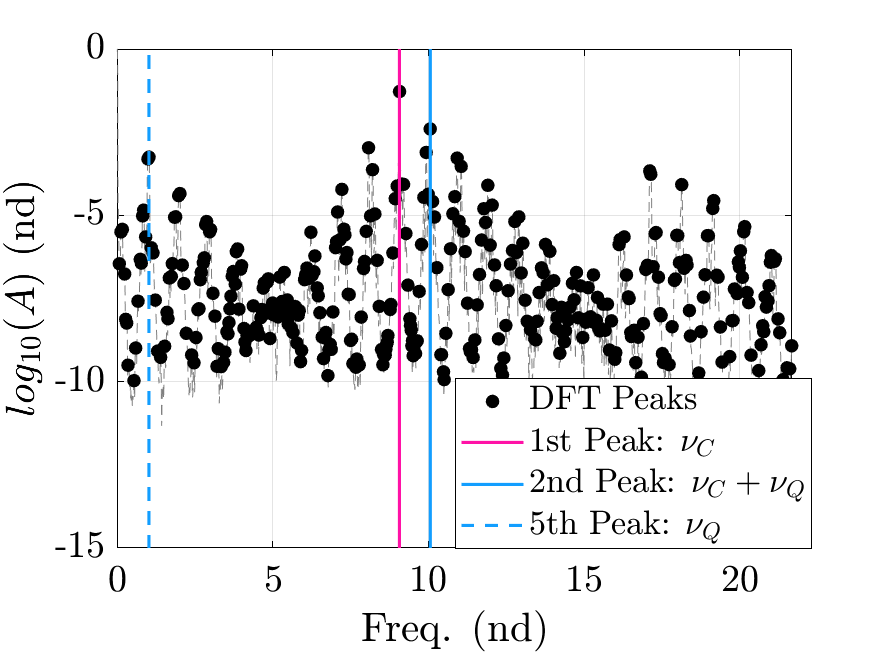}
        \caption{\acrshort{dft} and refined peaks with $\signalt = x$ (\acrshort{laskar})}
        \label{fig:dro-hfem-dft}
    \end{subfigure}
  \caption{\acrshort{hfem} quasi-\acrshort{dro} initiated with the state in Table \ref{table:dro-ig}}
  \label{fig:dro-hfem}
\end{figure}

The \acrfull{dro} generally displays stable behavior across dynamical models and serves as a good example to test the single-shooting formulation. One such example is illustrated in Figs. \ref{fig:dro-cr3bp-geometry} and \ref{fig:dro-hfem-geometry}, initiated with the state in Table \ref{table:dro_initial}. The initial state is propated for 20 years in the \acrshort{cr3bp} (Fig. \ref{fig:dro-cr3bp-geometry}) and the \acrshort{hfem} (Fig. \ref{fig:dro-hfem-geometry}). Note that Eq. \eqref{eq:mci_brf} is leveraged to rotate the \acrshort{brf} state into the \acrshort{mci}. An initial epoch of $2460576.5$ in Julian Date (09/23/2024) is leveraged. Both \acrshort{cr3bp} and \acrshort{hfem} display quasi-periodic behaviors, consistent with a previous observation by Bezrouk and Parker \cite{bezrouk2017long}. Note that the distance from the Moon to the $x$-axis crossing is approximately $22000$ km. 

An arbitrary goal is introduced to test the functionality of \acrshort{fdc}. Note that while the \acrshort{cr3bp} dynamics allow $(\inter=3)$-dimensional QPOs for DRO-related structures, for simplicity, consider an in-plane motion where $z = z' = 0$. Then, two fundamental frequencies exist within the \acrshort{cr3bp} framework: (1) $\fund{C}$, associated with the longitudinal period of the underlying \acrshort{dro}, and (2) $\fund{Q}$, associated with the in-plane quasi-periodic motion. The initial guess in Table \ref{table:dro-ig} is generated by deliberately perturbing a state on a periodic \acrshort{dro}; the initial value $x' = 0.01$ does not satisfy the perpendicular crossing condition ($x' = 0$). The \acrshort{dft} on $x$-position components illustrate this behavior in Fig. \ref{fig:dro-cr3bp-dft}, sampled with $\sampleSize = 2^{16}$. In contrast to the \acrshort{dft} results for a \acrfull{po} (Fig. \ref{fig:dro-compare-l-g-dft.png}), the peaks now involve more than one fundamental frequency. Refining the first two peaks ($1 \leq \idxA \leq 2$) with the \acrshort{laskar} algorithm, it is possible to identify the fundamental frequencies. The first peak $\freqA{1}$ coincides with $\fund{C}$, and the second peak $\freqA{2}$ coincides with $\fund{Q}$. The values for these fundamental frequencies are also confirmed from the monodromy matrix of the \acrshort{po} associated with the frequency of $\fund{C}$. The in-plane center eigenvalue is $\lambda_Q \approx 0.7362 \pm \imag 0.6767$ associated with a rotation number of $\sigma_Q = \arctan 2 (| \text{Im} (\lambda_Q)|, \text{Re} (\lambda_Q)) \approx 0.7433$ (rad). Noting that the rotation number is connected to the frequency as $\sigma_Q = 2\pi\fund{Q}/\fund{C}$, the second fundamental frequency from the monodromy matrix is retrieved as $\fund{Q} \approx 1.0249 \approx \freqA{2}$ (Fig. \ref{fig:dro-cr3bp-dft}). Of course, depending on the magnitude along the quasi-periodic oscillation, $\fund{Q}$ may coincide with the second dominant peak. Such a case is often observed within the \acrshort{hfem}. Figure \ref{fig:dro-hfem} illustrates the counterpart behavior within the \acrshort{hfem}. The initial \acrshort{brf} state from Table \ref{table:dro_initial} is rotated into the \acrshort{mci} on 09/23/2024, propagated within the \acrshort{hfem} and rotated back to \acrshort{brf} for consistent visualization. With the external perturbations, the \acrshort{dro} now tracks a more dispersed behavior as evidenced in Fig. \ref{fig:dro-hfem-geometry}. The corresponding Fourier domain response is recorded in Fig. \ref{fig:dro-hfem-dft} with $\signalt = x$ as consistently defined with the \acrshort{cr3bp} investigation from Figure \ref{fig:dro-cr3bp-dft}. The most dominant peak, $\freqA{1}$, still coincides with $\fund{C}$. Note that the exact value for $\fund{C}$ may be subject to change from the \acrshort{cr3bp} value as it is not explicitly controlled in the rotation process (Eq. \eqref{eq:mci_brf}. In contrast to the \acrshort{cr3bp} case, the \acrshort{hfem} trajectory does not supply $\freqA{2} \approx \fund{Q}$. Rather, $\freqA{2}$ is associated with the convoluted frequency between $\fund{C}$ and $\fund{Q}$. Refining frequencies associated with subsequent peaks, $\fund{Q}$ now coincides with $\freqA{5}$, the fifth dominant frequency. While it is more trivial to locate and target the most dominant frequency, $\fund{C}$ in this case, it is more nuanced to detect the secondary frequency, i.e., $\fund{Q}$. The situation is more challenging within the \acrshort{hfem} due to the presence of the external frequencies that generate multiple peaks that do not originally exist within the \acrshort{cr3bp} dynamics. The frequency structures in Figs. \ref{fig:dro-cr3bp-dft} and \ref{fig:dro-hfem-dft} are summarized in Tables \ref{table:cr3bp-dro-freq-initial} and \ref{table:hfem-dro-freq-initial}, respectively. 

Assume that the fundamental frequencies, $\fund{C}$ and $\fund{Q}$, are desired to be adjusted. The following (arbitrary) frequency structure is enforced, 
\begin{align}
    \label{eq:dro-target-structure} \text{For } \fund{C} : \quad \mb \freqA{\idxA} - 8.66331279836 \\ \phaseA{\idxA} - \pi \me = 0, \quad \text{For } \fund{Q} : \quad \mb \ampA{\idxA} - 0.001 \\ \phaseA{\idxA} - \pi/8 \me = 0.
\end{align}
Within the \acrshort{cr3bp}, examining the spectrum from Fig. \ref{fig:dro-cr3bp-dft}, $\freqA{\idxA}$ with $\idxA = 1, 2$ are \textit{identified} as $\fund{C}$ and $\fund{Q}$, respectively. During the corrections process within \acrshort{fdc}, however, the ordering of the secondary peaks may change. Alternatively, thus, the index $\idxA$ is recovered corresponding to the frequency bin that displays the $\idxA$-th peak, i.e., $\freqP{\idxA}$, closest to $\fund{Q} \approx 1.02597$ (\acrshort{hfem}). As the value for $\fund{Q}$, of course, changes during the corrections, these two strategies, i.e., (1) fixing the index $\idxA$ and (2) investigating the frequency $\fund{Q}$ supplied \textit{a priori}, are complementary to each other. Also note that the convoluted frequency $\fund{C} + \fund{Q}$ may be dominant compared to the fundamental frequency $\fund{Q}$ itself as apparent within the \acrshort{hfem} (Fig. \ref{fig:dro-hfem-dft}). It is possible to employ \acrshort{fdc} on convoluted frequencies, an approach not explored in the current investigation.

\begin{figure}[h!]
    \centering
    \begin{subfigure}[b]{0.48\textwidth}
        \centering
        \includegraphics[width=0.99\linewidth]{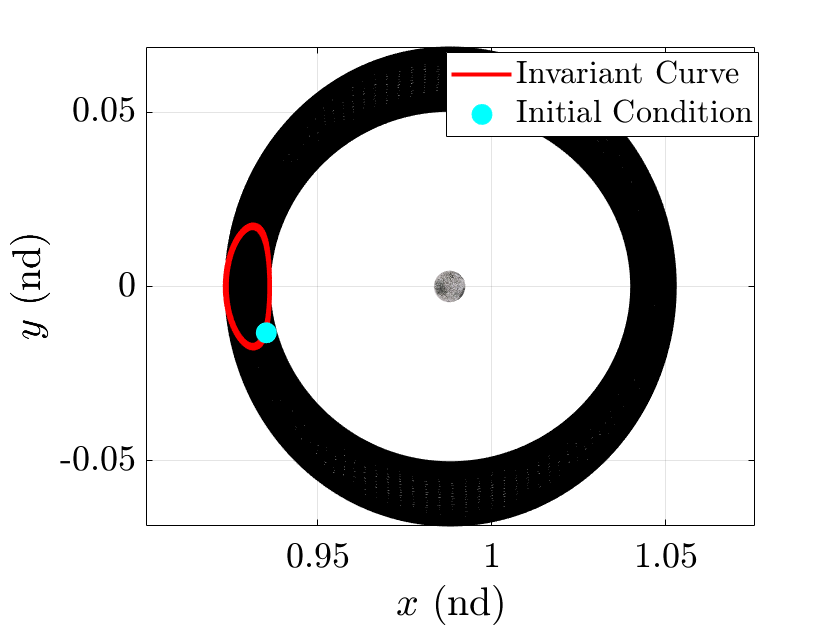}
        \caption{\acrshort{cr3bp} (\acrshort{brf})}
        \label{fig:dro-cr3bp-converged-geometry}
    \end{subfigure}
    \begin{subfigure}[b]{0.48\textwidth}
        \centering
        \includegraphics[width=0.99\linewidth]{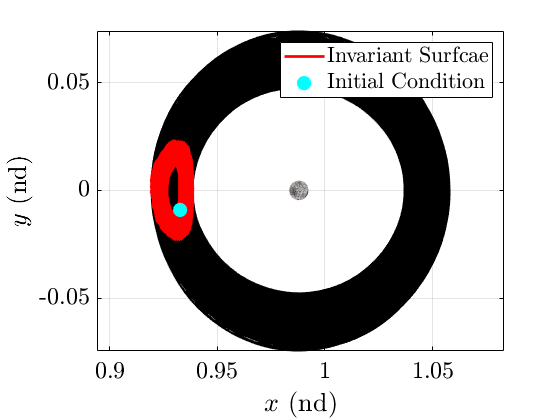}
        \caption{\acrshort{hfem} (\acrshort{brf})}
        \label{fig:dro-hfem-converged-geometry}
    \end{subfigure}
    \caption{Targeted quasi-\acrshort{dro}}\label{fig:dro-converged-geometry}
\end{figure}

The \acrshort{fdc} process is employed for both the \acrshort{cr3bp} and \acrshort{hfem} to satisfy the target frequency structure from Eq. \eqref{eq:dro-target-structure}. The resulting geometries are illustrated in Fig. \ref{fig:dro-converged-geometry}. The red curves within the plots correspond to the invariant curve (\acrshort{cr3bp}) and invariant surface (\acrshort{hfem}), respectively, recorded every passage of $2\pi/\fund{C}$ in the nd time, eliminating the dependency on one of the fundamental frequencies. As the phase angle for $\fund{C}$ is targeted as $\phaseA{\idxA} = \pi$, the invariant curve (surface) appears to the left of the Moon (Eq. \eqref{eq:signal}). The corresponding frequency structures are included in Tables \ref{table:cr3bp-dro-freq-target} and \ref{table:hfem-dro-freq-target}. The targeted frequency components are marked with \yes, converged below the set tolerance, $1\cdot 10^{-10}$ (nd). The initial state for the \acrshort{cr3bp} and \acrshort{hfem} along the targeted trajectories are included in Tables \ref{table:dro-cr3bp-converged} and \ref{table:dro-hfem-converged}, respectively. Compared to the original geometries from Figs. \ref{fig:dro-cr3bp-geometry} and \ref{fig:dro-hfem-geometry}, the targeted frequency structures are associated with larger $\ampA{\idxA}$ for $\fund{Q}$ and appear to be more dispersed within the \acrshort{brf}. The current example on the \acrshort{dro} demonstrates the capabilities of \acrshort{fdc} to specify trajectories on a QPO across models of varying fidelity. Specifically, the \acrshort{hfem} trajectory likely tracks higher-dimensional QPO with $\ttl \geq 7$. The proposed algorithm, i.e., \acrshort{fdc}, is capable of targeting $\inter$ (internal) frequencies, avoiding the complexity of explicitly considering $\ttl$ frequencies simultaneously. Note that the out-of-plane quasi-periodic motion also exists that is not explicitly considered in this example. Without an explicit control of such components, the corrections process in \acrshort{fdc} may introduce components into the third direction. Note that while $\x{0}$, the free variable vector, is six-dimensional, $\f{f}$ (Eq. \eqref{eq:frequency_constraint}) is a four-dimensional vector. A minimum-norm approach is leveraged to supply the search direction for the differential corrector.

\begin{table}[h!]
\centering
\caption{\label{table:dro-ig}Initial state $\x{0}$ for the quasi-\acrshort{dro} (Figs. \ref{fig:dro-cr3bp} and \ref{fig:dro-hfem})}
\begin{tabular}{|l|l|}
\hline
Component & Value \\ \hline
$x$ &  $0.929817046666844$     \\ \hline
$y$ &  $0$     \\ \hline
$z$ &  $0$     \\ \hline
$x'$ &  $0.01$      \\ \hline
$y'$ &  $0.522717065584611$      \\ \hline
$z'$ &  $0$     \\ \hline
\end{tabular}
\end{table}

\begin{table}[h!]
\centering
\caption{\label{table:cr3bp-dro-freq-initial}Initial frequency structure for the \acrshort{cr3bp} quasi-\acrshort{dro} propagated with Table \ref{table:dro-ig}}
\begin{tabular}{|l|l|l|}
\hline
 & $\fund{C}$ ($\idxA = 1$) & $\fund{Q}$ ($\idxA = 2$)  \\ \hline
$\freqA{\idxA}$ & $8.663312798420872$ &$1.024930860632975$  \\ \hline
$\ampA{\idxA}$ & $0.058190266942043$ & $0.001425616315091$ \\ \hline
$\phaseA{\idxA}$ & $-3.098796017720698$ & $1.571460016637346$ \\ \hline
\end{tabular}
\end{table}

\begin{table}[h!]
\centering
\caption{\label{table:hfem-dro-freq-initial}Initial frequency structure for the \acrshort{hfem} quasi-\acrshort{dro} propagated with Table \ref{table:dro-ig}}
\begin{tabular}{|l|l|l|}
\hline
 & $\fund{C}$ ($\idxA = 1$)& $\fund{Q}$ ($\idxA = 5$) \\ \hline
$\freqAArg$ & $9.067270092334322$ &$1.025971597715591$  \\ \hline
$\ampA{}$ & $0.056230858215909$ & $0.000563726693198$ \\ \hline
$\phaseA{}$ & $3.062426060730863$ & $2.257364738089462$ \\ \hline
\end{tabular}
\end{table}

\begin{table}[h!]
\centering
\caption{\label{table:cr3bp-dro-freq-target}Targeted frequency structure for the \acrshort{cr3bp} quasi-\acrshort{dro} propagated with Table \ref{table:dro-cr3bp-converged} (\yes: targeted)}
\begin{tabular}{|l|l|l|}
\hline
 & $\fund{C}$ ($\idxA = 1$)& $\fund{Q}$ ($\idxA = 2$)  \\ \hline
$\freqA{\idxA}$ &\textbf{${8.663312798369873}$} \yes &$1.025840351634718$  \\ \hline
$\ampA{\idxA}$ & $0.057630566792085$ & \textbf{${0.009999999891364}$} \yes \\ \hline
$\phaseA{\idxA}$ & \textbf{${3.141592655033092}$} \yes & \textbf{${0.785398163502869}$} \yes \\ \hline
\end{tabular}
\end{table}

\begin{table}[h!]
\centering
\caption{\label{table:hfem-dro-freq-target}Targeted frequency structure for the \acrshort{hfem}  quasi-\acrshort{dro} propagated with Table \ref{table:dro-hfem-converged} (\yes: targeted)}
\begin{tabular}{|l|l|l|}
\hline
 & $\fund{C}$ ($\idxA = 1$)& $\fund{Q}$ ($\idxA = 5$)  \\ \hline
$\freqA{\idxA}$ & \textbf{${8.663312798372578}$} \yes &$1.026122000070913$  \\ \hline
$\ampA{\idxA}$ & $0.057532681072891$ & \textbf{${0.009999999998136}$} \yes \\ \hline
$\phaseA{\idxA}$ & \textbf{${3.141592653498086}$} \yes & \textbf{${0.785398163697401}$} \yes \\ \hline
\end{tabular}
\end{table}

\begin{table}[h!]
\centering
\caption{\label{table:dro-cr3bp-converged}Targeted initial state $\x{0}$ for the quasi-\acrshort{dro}, \acrshort{cr3bp} (Fig. \ref{fig:dro-cr3bp-geometry})}
\begin{tabular}{|l|l|}
\hline
Component & Value \\ \hline
$x$ &  $0.935182071235081$     \\ \hline
$y$ & $-0.013296687505060$ \\ \hline
$z$ &  $0$     \\ \hline
$x'$ &  $ -0.071999754120518 $      \\ \hline
$y'$ &  $0.542602051432254$      \\ \hline
$z'$ &  $0$     \\ \hline
\end{tabular}
\end{table}

\begin{table}[h!]
\centering
\caption{\label{table:dro-hfem-converged}Targeted initial state $\x{0}$ for the quasi-\acrshort{dro}, \acrshort{hfem} (Fig. \ref{fig:dro-hfem-geometry})}
\begin{tabular}{|l|l|}
\hline
Component & Value \\ \hline
$x$ &  $0.932714478796438$     \\ \hline
$y$ &  $-0.009071158291536 $     \\ \hline
$z$ &  $-0.000720383507459$     \\ \hline
$x'$ &  $-0.018764477779327$      \\ \hline
$y'$ &  $0.543649830698251$      \\ \hline
$z'$ &  $0.000579106771136$     \\ \hline
\end{tabular}
\end{table}

\pagebreak

\subsection{\acrfull{elfo} Constellation Design}

\acrfull{elfos} offer an attractive option to support various missions in the lunar vicinity. These orbital structures supply near-stability over a long horizon time and the ability to provide extended line-of-sight to the Lunar South Pole (LSP), a key region for expanding cislunar space activities in the near future \cite{baker2024comprehensive}. Leveraging these benefits, multiple space entities envision the inclusion of one or more satellites in \acrshort{elfos}\footnote{\url{https://www.intuitivemachines.com/post/nasa-selects-intuitive-machines-to-deliver-4-lunar-payloads-in-2024}, Accessed: 09-05-2024}\footnote{\url{https://spacenews.com/china-could-develop-dual-relay-satellite-system-for-earth-moon-communications-to-reduce-geopolitical-risks}, Accessed: 09-05-2024} as a part of a data relay system. The \acrshort{elfos} satisfy the ``frozen'' condition, admitting nearly constant \textit{mean} eccentricity, inclination, and argument of perilune \cite{nie2018lunar}. While a variety of orbits supply such a condition, the \acrshort{elfos} are characterized by relatively high altitude, inclination and eccentricity values, where the point-mass Earth gravity significantly contributes to the stable evolution of the orbit geometry \cite{ely2005stable}. 

The \acrshort{elfos} exist in the vicinity of the equilibria derived from the simplified, \acrfull{dam} \cite{nie2018lunar, ely2005stable, folta2006lunar, de2003third2}. From the equations of motion as represented in the doubly-averaged orbital elements from Eqs. \eqref{eq:dadt}-\eqref{eq:dthetadt}, observe that the dynamics for $\da{e}, \da{i}, \da{\omega}$ are separable. Then, the following equilibrium conditions specify the \acrshort{elfos},
\begin{align}
    \label{eq:omega_eq} &\da{\omega}_{eq} = \pi/2, 3\pi/2  \\
    \label{eq:ecc_eq} &\da{e}_{eq} = \sqrt{1-5/3 \cos^2 \da{i}_{eq} }.
\end{align}
The equilibria are elliptic, i.e., they admit nearby periodic orbits in terms of $\da{e}, \da{i}, \da{\omega}$. This period is denoted as the \textit{long-period} and constitutes one of the fundamental frequencies, $\fund{L}$, for the \acrshort{elfos}. While specific values differ by the combinations of $\da{e}-\da{i}-\da{\omega}$, it is generally associated with periods on the order of one year. Recall that the \acrshort{dam} is integrable in the Liouville-Arnold sense \cite{arnol2013mathematical, nie2018lunar} where every solution generally appears as an $(\inter=3)$-dimensional quasi-periodic orbit. Thus, two more fundamental frequencies exist that are associated with the averaged angles. The \textit{short-period} corresponds to one complete revolution for $\da{M}$, that depends on $\da{a}$, and is on the order of one day. The \textit{medium-period} is associated with the passage of $2\pi$ (radians) for the modified right ascension defined as $\da{\Omega}_M = t - \da{\Omega}$ as measured clockwise within the \acrshort{brf} such that the rate of $\da{\Omega}_M$ is positive. Noting that $d \da{\Omega} / dt < 0$ for $\da{\omega} = \pi/2$ and $\da{i} \leq \pi/2$ (Eq. \eqref{eq:dOmegadt}), the period for $\da{\Omega}_M$ is slightly larger than $2\pi$, on the order of one month. The fundamental frequencies for the \acrshort{elfos} and structures in the near vicinity are, thus, comprised of $\fund{S}, \fund{M}, \fund{L}$, associated with the short-, medium-, and long-period of the system. The \acrshort{elfos} pertain to the specific case where the Eqs. \eqref{eq:omega_eq} and \eqref{eq:ecc_eq} are satisfied and the long-period frequency $\fund{L}$ is defined only in a limiting sense. 

While the \acrshort{elfos} are initially derived in the simplified \acrshort{dam}, the trajectories must be \textit{refined} within the \acrfull{hfem} such that they are suitable for actual flight. One challenge exists in defining a suitable transformation between the doubly-averaged orbital elements and the osculating elements (Fig. \ref{fig:frame_transformation}). While both semi-analytical \cite{nie2018lunar} and fully-numerical strategies \cite{howell2007design, ely2006constellations} appear in the literature to supply the transformation, a systematic corrections process remains elusive to refine \acrshort{elfos} within a higher-fidelity model. For designing and operating a satellite \textit{constellation} within \acrshort{elfos}, the capabilities to refine \acrshort{elfos} within a higher-fidelity model are increasingly important. The constellation introduces phase constraints between the satellites, necessitating precise initial states to maintain these constraints throughout the mission timeline.

\begin{figure}[h!]
    \centering
    \begin{subfigure}[b]{0.48\textwidth}
        \centering
        \includegraphics[width = 0.99\textwidth]{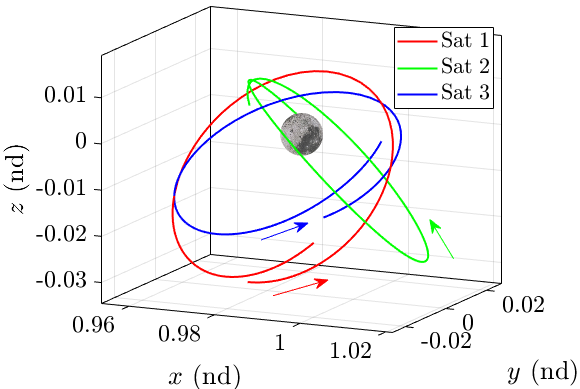}        \caption{\label{fig:three_sat_ig_config}Trajectories for the first $\approx1.3$ days (\acrshort{brf})}
    \end{subfigure}   
    \begin{subfigure}[b]{0.48\textwidth}
        \centering
        \includegraphics[width = 0.99\textwidth]{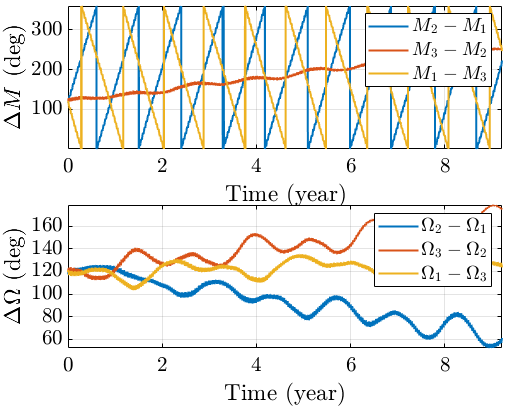}        \caption{\label{fig:three_sat_ig_drift}Phase drift in $M, \Omega$\protect\footnotemark} 
    \end{subfigure}
    \caption{\label{fig:three_sat_ig}Initial guess for the three-satellite constellation with undesired phase drift}
\end{figure}

The proposed \acrshort{fdc} approach is well-suited for a satellite constellation as the phasing information is explicitly enforced during the corrections process. The following numerical example demonstrates its capabilities within the context of the \acrshort{elfo} constellation. Consider an arbitrary constellation comprised of three satellites. The initial states are derived initially from the \acrshort{dam}, included in Table \ref{table:initial_condition}, designed to provide coverage of the lunar south pole ($\da{\omega} = 90^\circ$) for approximately 10 years. The combination of $\da{e},\da{i},\da{\omega}$ satisfy the conditions from Eqs. \eqref{eq:omega_eq} and \eqref{eq:ecc_eq}. Phasing requirements are set such that the satellites are separated by $120^\circ$ in both $\da{\Omega}$ and $\da{M}$, ensuring consistent coverage. While this requirement suggests initial guesses for $\da{\text{\oe}}$ (a symbol defined in Section \ref{sec:dynamics}) offset by $120^\circ$ in $\da{\Omega}$ and $\da{M}$, the transformed states $\vec{\mathbb{X}}_0$ (represented in the \acrshort{mci}) in the \acrshort{hfem} often fail to meet the desired phasing upon test propagation, requiring adjustments. For the examples presented, the initial epoch is fixed on 09/23/2024. The resulting trajectories are visualized in Fig. \ref{fig:three_sat_ig_config} leveraging the \acrshort{brf} as the reference frame. The plots depict the first $\approx 1.3$ days of propagation. Over longer periods, however, phase drift in $M$ and $\Omega$ is observed, as apparent in Fig. \ref{fig:three_sat_ig_drift}. The vertical axis illustrates $\Delta M$ and $\Delta \Omega$, that is, the relative phasing between the satellites within the constellation over time. This phase drift negatively impacts the constellation’s ability to maintain consistent coverage and line-of-sight.

\begin{figure}[h!]
    \centering
    \begin{subfigure}[b]{0.48\textwidth}
        \centering
        \includegraphics[width = 0.99\textwidth]{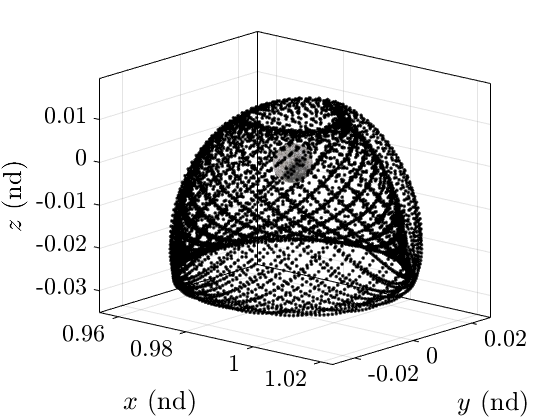}
        \caption{Geometry within the \acrshort{brf}, first $\approx 70$ days}
        \label{fig:elfo-hfem-geometry}
    \end{subfigure} \\
    \begin{subfigure}[b]{0.48\textwidth}
        \centering
        \includegraphics[width = 0.99\textwidth]{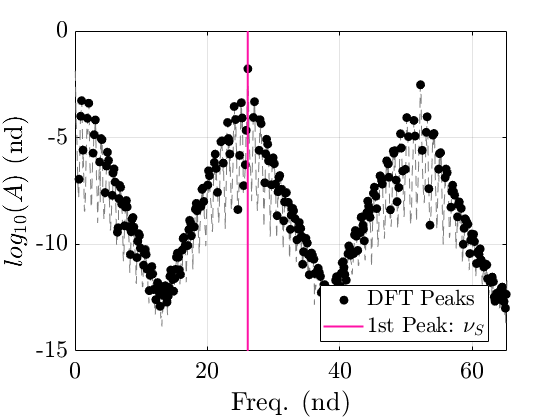}
        \caption{\acrshort{dft} and refined peaks with $\signalt = z$ (\acrshort{gomez})}
        \label{fig:elfo-hfem-dft-z}
    \end{subfigure}
    \begin{subfigure}[b]{0.48\textwidth}
        \centering
        \includegraphics[width = 0.99\textwidth]{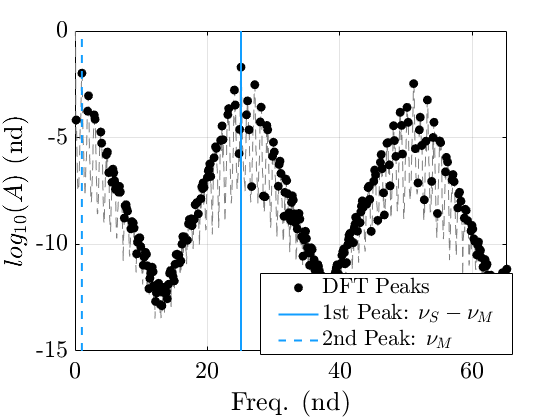}
        \caption{\acrshort{dft} and refined peaks with $\signalt = x$ (\acrshort{gomez})}
        \label{fig:elfo-hfem-dft-x}
    \end{subfigure}
  \caption{\acrshort{hfem} \acrshort{elfo} (Sat 1) initiated with the state in Table \ref{table:contellation_numbers}}
  \label{fig:elfo-hfem}
\end{figure}

\begin{figure}[h!]
    \centering
    \begin{subfigure}[b]{0.48\textwidth}
        \centering
        \includegraphics[width = 0.99\textwidth]{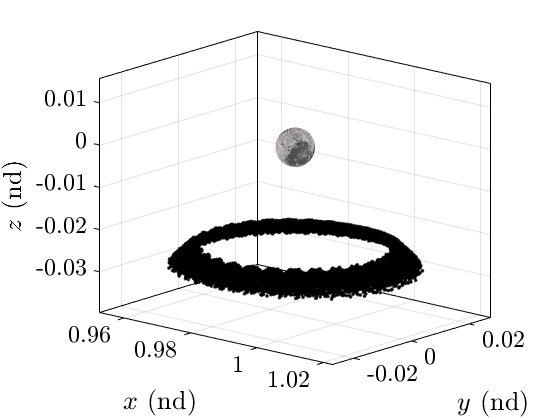}
        \caption{Geometry within the \acrshort{brf}, first $2^{14}$ points}
        \label{fig:elfo-hfem-geometry-strobo}
    \end{subfigure} 
    \begin{subfigure}[b]{0.48\textwidth}
        \centering
        \includegraphics[width = 0.99\textwidth]{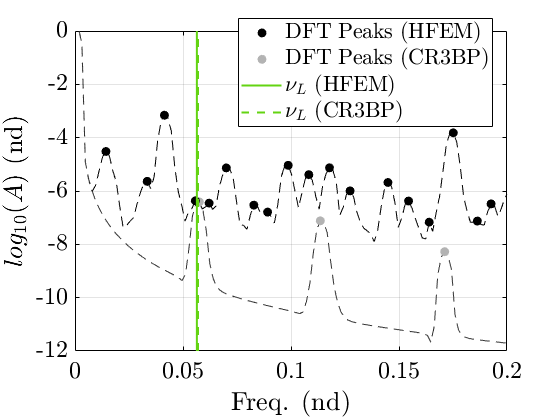}
        \caption{\acrshort{dft} and refined peaks from $\signalt = x$ (\acrshort{gomez})}
        \label{fig:elfo-hfem-dft-x-strobo}
    \end{subfigure}
  \caption{Stroboscopic map along the \acrshort{hfem} \acrshort{elfo} (Sat 1) initiated with the state in Table \ref{table:contellation_numbers}}
  \label{fig:elfo-hfem-strobo}
\end{figure}

The first satellite, i.e., Sat 1, is further examined to supply insights regarding the geometric interpretations for the drift and the fundamental frequencies, i.e., $\fund{S}, \fund{M}, \fund{L}$. The geometry is depicted within \acrshort{brf} in Fig. \ref{fig:elfo-hfem-geometry} for the first $\approx10$ years. Note that the orbital plane rotates clockwise within the \acrshort{brf}. As the short-period frequency $\fund{S}$ is closely related to $\da{M}$, and from the shape of the orbit, $\signalt = z$ is examined. The resulting spectrum is plotted in Fig. \ref{fig:elfo-hfem-dft-z}, where the most dominant peak tracks $\fund{S} \approx 26$ (nd), refined with the \acrshort{gomez} method. However, as the medium-period frequency, $\fund{M}$, is associated with $\da{\Omega}_M$, it is rather challenging to extract the information directly from $z$. Thus, $\signalt = x$ is examined and plotted in Fig. \ref{fig:elfo-hfem-dft-x}, where the second peak ($\freqA{\idxA}$) tracks $\fund{M}$, also refined with the \acrshort{gomez} algorithm. The initial condition derived within the \acrshort{dam} satisfies the equilibrium condition from Eqs. \eqref{eq:omega_eq}-\eqref{eq:ecc_eq}, eliminating the peaks associated with the long-period frequency $\fund{L}$ within the \acrshort{dam} (tracking a 2-dimensional \acrshort{qpo}). However, due to the limited accuracy of the transformation process described in Fig. \ref{fig:frame_transformation}, trajectories within the \acrshort{hfem} typically excite the oscillatory mode into the $\fund{L}$-direction. This long-period (low) frequency is associated with a low amplitude, requiring extra caution in proper detection. One requirement is a long propagation time such that $\fndS \ll 1$ to resolve the low frequency $\fund{L} < 0.1$ (on the order of years). Note that this lower-frequency-domain is challenging to visually detect within the spectra included in Fig. \ref{fig:elfo-hfem}. One possible complementary strategy is to leverage a stroboscopic map with a period of $2\pi/\fund{S}$, eliminating the influence of the short-period frequency $\fund{S}$. The states at the stroboscopic map are illustrated in Fig. \ref{fig:elfo-hfem-geometry-strobo}, marking the first $\sampleSize = 2^{14}$ returns to the map. The \acrshort{dft} results are included in Fig. \ref{fig:elfo-hfem-dft-x-strobo} with $\signalt = x$. While the \acrshort{hfem} alone demonstrates multiple peaks due to the impact from the $\ex$ external frequencies, the auxiliary analysis from the \acrshort{cr3bp} initiated with the same initial state (from Table \ref{table:contellation_numbers}, transformed with Fig. \ref{fig:frame_transformation}) supplies well-defined peaks and a reference number for $\fund{L}\approx 0.057$ (nd).

\begin{table}[htpb]
\centering
\caption{\label{table:initial_condition}Set of doubly-averaged mean orbital elements for the three-satellite constellation}
\begin{tabular}{|c|c|c|c|c|c|c|}

\hline
$\da{\text{\oe}}$ & $\da{a}$ & $\da{e}$  & $\da{i}$ & $\da{\Omega}$ & $\da{\omega}$ & $\da{M}$ \\ \hline
Value & $10,000$ (km)  & $0.4082$ (nd)  & $45^\circ$ & $0^\circ, 120^\circ, 240^\circ$ &  $90^\circ$ & $180^\circ, 300^\circ, 60^\circ$ \\ \hline
\end{tabular}
\end{table}

\begin{figure}[h!]
    \centering
    \begin{subfigure}[b]{0.48\textwidth}
        \centering
        \includegraphics[width = 0.99\textwidth]{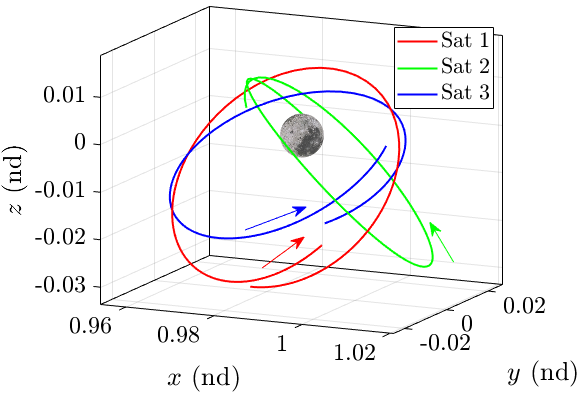}        \caption{\label{fig:three_sat_converged_config}Trajectories for the first $\approx1.3$ days (\acrshort{brf})}
    \end{subfigure}   
    \begin{subfigure}[b]{0.48\textwidth}
        \centering
        \includegraphics[width = 0.99\textwidth]{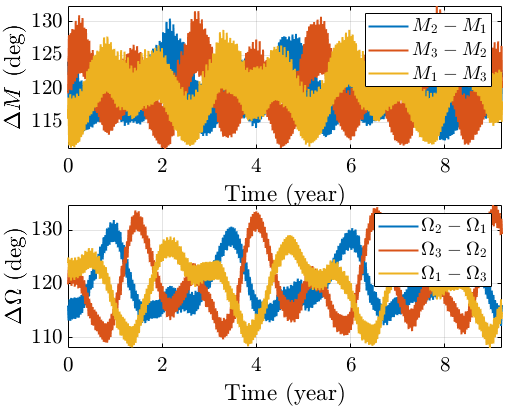}        \caption{\label{fig:three_sat_converged_drift}Phase drift in $M, \Omega$}
    \end{subfigure}  
    \caption{\label{fig:three_sat_converged}Targeted and converged three-satellite \acrshort{elfo} constellation with controlled phase drift}
\end{figure}

\footnotetext{Note that the desired phase drift values for $M, \Omega$ are $120^\circ$ on average. The subscripts $1,2,3$ denote Sat 1, 2, 3; for example, $M_2 - M_1$ corresponds to the phase difference in $M$ between Sat 2 with respect to Sat 1. The same description applies to Fig. \ref{fig:three_sat_converged_drift} as well. }

\clearpage

With this geometric interpretation of the fundamental frequencies, the initial configuration for the \acrshort{elfo} constellation from Fig. \ref{fig:three_sat_ig} is revisited. The initial state and frequency structure for each of the satellites are included in Table \ref{table:contellation_numbers}. Note that the initial state $\x{0}$ is represented in the inertial frame (\acrshort{mci}), but the signals for the \acrshort{dft} are constructed within the rotating frame (\acrshort{brf}), e.g., $x$. Thus, the partial in Eq. \eqref{eq:dsignal_dx} is evaluated with,
\begin{align}
    & \signalt = x = \mb 1, 0, 0\me \vec{\rho} = \mb 1, 0, 0\me \frac{l_*}{l} \bm{C}^T\vec{R} + \vec{\rho}_M,
\end{align}
from Eq. \eqref{eq:mci_brf}. Then, 
\begin{align}
    & \frac{d \signalt(\tidxT)}{d \x{0}} =  \frac{d \signalt(\tidxT)}{d \xp{0}{\tidxT}{\x{0}}{\epoch{0}}} \frac{d \xp{0}{\tidxT}{\x{0}}{\epoch{0}}}{d \x{0}} = \mb 1, 0, 0\me \frac{l_*}{l} \bm{C}^T \stm{0}{\tidxT}{\x{0}}{\epoch{0}}_{3\times 6} .
\end{align}
The subscript ${3\times6}$ denotes that a sub-matrix of the STM is leveraged that is connected to the position at sample time $\tidxT$. The direction cosine matrix $\bm{C}$ is also evaluated at time $\tidxT$. A similar process is extended to $\signalt = z$ as well. The frequency structures included in Table \ref{table:contellation_numbers} consist of five different quantities: frequency and phase angles for short- and medium-period frequencies, and the amplitude associated with the long-period frequency. Note that these fundamental frequencies and the associated phase angles are retrieved from the frequency refinement analysis in Fig. \ref{fig:elfo-hfem}. Arrows indicate whether components overshoot ($\lag$) or undershoot ($\lead$) their desired values with respect to Sat 1. For example, Sat 1 employs an initial guess $\nu_S = 26.0840$ (nd), while Sat 2 and Sat 3 exhibit slightly higher values (overshot, $\lag$). The initial design intends for Sat 2 and Sat 3 to lead and lag Sat 1 by $120^\circ$, respectively. However, these small deviations in $\nu_S$ prevent these satellites from maintaining the desired phase offset, as evident in Fig. \ref{fig:three_sat_ig_drift}. Similar behavior is observed for $\nu_M$, where the relative configuration for ${\Omega}$ deteriorates over time. Note that $\theta_S \approx \da{M}$ and $\theta_M \approx \da{\Omega}_M$,  although not exactly. Furthermore, the initial phase offsets for $\nu_S$ and $\nu_M$ deviate slightly from the desired $120^\circ$ due to limitations in the transformation accuracy described in Fig. \ref{fig:frame_transformation}. For instance, Sat 2 initially leads Sat 1 in $\theta_S$ by only $106^\circ$, falling short by $14^\circ$ for the intended offset (undershot $\lead$). For all three satellites, the trajectories are associated with small amplitudes for the long-period frequency $\fund{L}$ examined with $\signalt = x$, i.e., $\ampA{L} <1\cdot 10^{-6}$ (nd). The phase angles for this oscillatory mode are not illustrated in Table \ref{table:contellation_numbers}. 

Utilizing the frequency-domain targeter, the state vectors $\vec{\mathbb{X}}_0$ are refined to eliminate the phase drift. The following constraint vector is arbitrarily enforced,
\begin{align}
    \label{eq:F_min_norm}\text{For }\fund{S}:\quad \mb \freqA{\idxA} - 26.1685 \\ \Delta \phaseA{\idxA} -( \pm 120^\circ) \me, \quad \text{For }\fund{M}:\quad \mb \freqA{\idxA} - 1.0341 \\ \Delta \phaseA{\idxA} - (\pm 120^\circ)  \me, \quad \text{For }\fund{L}:\quad \mb \ampA{\idxA} - 1\cdot 10^{-6} \me <0 .
\end{align}
For Sat 1, the phase angles $\theta_{S, M}$ are not targeted as it is assumed to serve as an arbitrary reference satellite within the constellation. After Sat 1 is targeted with Eq. \eqref{eq:F_min_norm}, the remaining satellites are adjusted relative to Sat 1 with phase offsets of $\pm 120^\circ$. The long-period frequency component, $\nu_L$, is not explicitly targeted. Rather, it is assumed that the quasi-periodic oscillation amplitude in $x$ associated with $\nu_L \approx 0.057$ (nd) is desired to be capped at $1 \cdot 10^{-6}$ (nd). Solving Eq. \eqref{eq:F_min_norm} using a minimum-norm approach ensures that $\ampA{\idxA}$ (the oscillation magnitude in $x$) does not significantly deviate from its original value for $\fund{L}$; while $\ampA{\idxA}$ is monitored for $\fund{L}$, it is not explicitly constrained in this example. The resulting state vectors $\vec{\mathbb{X}}_0$ are provided in Table \ref{table:contellation_numbers}. The frequency structures are now fully targeted successfully, as indicated by green check marks (\yes). The targeted frequencies $\nu_S$ and $\nu_M$ match within a tolerance of $1 \cdot 10^{-7}$ (nd), and the initial phase angles are set to be $120^\circ$ apart with a tolerance of $1 \cdot 10^{-7}$ radians. This information is confirmed by examining the evolution of the relative $M, \Omega$ angles over $\approx 10$ years as plotted in Fig. \ref{fig:three_sat_converged_drift}. While the oscillations exist due to the nature of the osculating elements (\oe), it is obvious that the relative phasing across the constellation remains stable, ensuring consistent performance.

\begin{table}[htpb]
\centering
\caption{\label{table:contellation_numbers}\acrshort{elfo} constellation constructed on 09/23/2024 (\yes: targeted, \lag: overshoot, \lead: undershoot) \tablefootnote{The angles $\Delta \theta_S$, $\Delta \theta_M$ are defined with respect to Sat 1.}}
\begin{tabular}{|c|ccc|ccc|}
\hline
                       & \multicolumn{3}{c|}{Initial guess}                                                                                                                                                                                    & \multicolumn{3}{c|}{Targeted and converged}                                                                                                                                                                                    \\ \hline
\multicolumn{1}{|c|}{Sat No.} & \multicolumn{1}{c|}{Sat 1}                                               & \multicolumn{1}{c|}{Sat 2}                                                     & \multicolumn{1}{c|}{Sat 3}                                & \multicolumn{1}{c|}{Sat 1}                                               & \multicolumn{1}{c|}{Sat 2}                                                   & \multicolumn{1}{c|}{Sat 3}                              \\ \thickhline
$\nu_S$ (nd)               & \multicolumn{1}{c|}{\begin{tabular}[c]{@{}c@{}}26.0840\\ (N/A)\end{tabular}} & \multicolumn{1}{c|}{\begin{tabular}[c]{@{}c@{}}$26.1687$\\ ($\lag 0.09$)\end{tabular}} & \begin{tabular}[c]{@{}c@{}}$26.1685$\\ ($\lag 0.09$)\end{tabular} & \multicolumn{1}{c|}{\begin{tabular}[c]{@{}c@{}}$26.1685$\\ (N/A)\end{tabular}} & \multicolumn{1}{c|}{\begin{tabular}[c]{@{}c@{}}$26.1685$\\ (\yes)\end{tabular}} & \begin{tabular}[c]{@{}c@{}}$26.1685$\\ (\yes) \end{tabular} \\ \hline
$\Delta \theta_S$ (deg)              & \multicolumn{1}{c|}{\begin{tabular}[c]{@{}c@{}}$0^\circ$\\ (N/A)\end{tabular}}                                                   & \multicolumn{1}{c|}{\begin{tabular}[c]{@{}c@{}}$106^\circ$\\ (\lead $14^\circ$)\end{tabular}}       & \begin{tabular}[c]{@{}c@{}}$242^\circ$\\ (\lag $2^\circ$)\end{tabular}        & \multicolumn{1}{c|}{\begin{tabular}[c]{@{}c@{}}$0^\circ$\\ (N/A)\end{tabular}}                                                  & \multicolumn{1}{c|}{\begin{tabular}[c]{@{}c@{}}$120^\circ$\\ (\yes)\end{tabular}}     & \begin{tabular}[c]{@{}c@{}}$240^\circ$\\ (\yes)\end{tabular}     \\ \hline
$\nu_M$ (nd)               & \multicolumn{1}{c|}{\begin{tabular}[c]{@{}c@{}}$1.0341$\\ (N/A)\end{tabular}}  & \multicolumn{1}{c|}{\begin{tabular}[c]{@{}c@{}}$1.0367$\\ ($\lag 0.026$)\end{tabular}} & \begin{tabular}[c]{@{}c@{}}$1.0327$\\ ($\lead0.014$)\end{tabular} & \multicolumn{1}{c|}{\begin{tabular}[c]{@{}c@{}}$1.0341$\\ (N/A) \end{tabular}}  & \multicolumn{1}{c|}{\begin{tabular}[c]{@{}c@{}}$1.0341$\\ (\yes)\end{tabular}}  & \begin{tabular}[c]{@{}c@{}}$1.0341$\\ (\yes)\end{tabular}  \\ \hline
$\Delta \theta_M$ (deg)              & \multicolumn{1}{c|}{\begin{tabular}[c]{@{}c@{}}$0^\circ$\\ (N/A)\end{tabular}}                                                    & \multicolumn{1}{c|}{\begin{tabular}[c]{@{}c@{}}$133^\circ$\\ ($\lag13^\circ$)\end{tabular}}       & \begin{tabular}[c]{@{}c@{}}$235^\circ$\\ ($\lead5^\circ$)\end{tabular}        & \multicolumn{1}{c|}{\begin{tabular}[c]{@{}c@{}}$0^\circ$\\ (N/A)\end{tabular}}                                                   & \multicolumn{1}{c|}{\begin{tabular}[c]{@{}c@{}}$120^\circ$\\ (\yes)\end{tabular}}     & \begin{tabular}[c]{@{}c@{}}$240^\circ$\\ (\yes)\end{tabular}     \\ \hline
$A_{L}$ (nd)               & \multicolumn{1}{c|}{$<1\cdot 10^{-6}$}                                                    & \multicolumn{1}{c|}{$<1\cdot 10^{-6}$}                                                          &    $<1\cdot 10^{-6}$                                                       & \multicolumn{1}{c|}{$<1\cdot 10^{-6}$}                                                    & \multicolumn{1}{c|}{$<1\cdot 10^{-6}$}                                                        &     $<1\cdot 10^{-6}$                                                    \\ \thickhline
$X$ (nd)               & \multicolumn{1}{c|}{\begin{tabular}[c]{@{}c@{}}$2.3681$\\ $\cdot 10^{-2}$ \end{tabular}}                                                    & \multicolumn{1}{c|}{\begin{tabular}[c]{@{}c@{}}$-2.4116$\\ $\cdot 10^{-2}$ \end{tabular}}                                         &   \begin{tabular}[c]{@{}c@{}}$-1.8681$\\ $\cdot 10^{-2}$ \end{tabular}                                                & \multicolumn{1}{c|}{\begin{tabular}[c]{@{}c@{}}$2.3046$\\ $\cdot 10^{-2}$ \end{tabular}}                                                    & \multicolumn{1}{c|}{\begin{tabular}[c]{@{}c@{}}$-2.4121$\\ $\cdot 10^{-2}$ \end{tabular}}                                                        &  \begin{tabular}[c]{@{}c@{}}$-1.8170$\\ $\cdot 10^{-2}$ \end{tabular}                                                       \\ \hline
$Y$ (nd)                & \multicolumn{1}{c|}{\begin{tabular}[c]{@{}c@{}}$3.9275$\\ $\cdot 10^{-3}$ \end{tabular}}                                                    & \multicolumn{1}{c|}{\begin{tabular}[c]{@{}c@{}}$4.3168$\\ $\cdot 10^{-3}$ \end{tabular}}                                                          &  \begin{tabular}[c]{@{}c@{}}$1.5963$\\ $\cdot 10^{-2}$ \end{tabular}                                                        & \multicolumn{1}{c|}{\begin{tabular}[c]{@{}c@{}}$4.5786$\\ $\cdot 10^{-3}$ \end{tabular}}                                                    & \multicolumn{1}{c|}{\begin{tabular}[c]{@{}c@{}}$4.4484$\\ $\cdot 10^{-3}$ \end{tabular} }                                                        &  \begin{tabular}[c]{@{}c@{}}$1.6671$\\ $\cdot 10^{-2}$ \end{tabular}                                        \\ \hline
$Z$ (nd)               & \multicolumn{1}{c|}{\begin{tabular}[c]{@{}c@{}}$-2.7631$\\ $\cdot 10^{-2}$ \end{tabular}}                                                    & \multicolumn{1}{c|}{\begin{tabular}[c]{@{}c@{}}$-3.4412$\\ $\cdot 10^{-3}$ \end{tabular} }                                                          &   \begin{tabular}[c]{@{}c@{}}$2.8722$\\ $\cdot 10^{-3}$ \end{tabular}                                                       & \multicolumn{1}{c|}{\begin{tabular}[c]{@{}c@{}}$-2.7283$\\ $\cdot 10^{-2}$ \end{tabular}}                                                    & \multicolumn{1}{c|}{\begin{tabular}[c]{@{}c@{}}$-2.4864$\\ $\cdot 10^{-3}$ \end{tabular}}                                                        &             \begin{tabular}[c]{@{}c@{}}$2.5967$\\ $\cdot 10^{-3}$ \end{tabular}                                \\ \hline
${X}'$ (nd)               & \multicolumn{1}{c|}{\begin{tabular}[c]{@{}c@{}}$1.6878$\\ $\cdot 10^{-1}$ \end{tabular} }                                                    & \multicolumn{1}{c|}{\begin{tabular}[c]{@{}c@{}}$1.4913$\\ $\cdot 10^{-1}$ \end{tabular}}                                                          &  \begin{tabular}[c]{@{}c@{}}$-4.5236$\\ $\cdot 10^{-1}$ \end{tabular}                                                          & \multicolumn{1}{c|}{\begin{tabular}[c]{@{}c@{}}$1.7087$\\ $\cdot 10^{-1}$ \end{tabular}}                                                    & \multicolumn{1}{c|}{\begin{tabular}[c]{@{}c@{}}$1.5159$\\ $\cdot 10^{-1}$ \end{tabular} }                                                        &  \begin{tabular}[c]{@{}c@{}}$-4.4563$\\ $\cdot 10^{-1}$ \end{tabular}                                           \\ \hline
$Y'$ (nd)               & \multicolumn{1}{c|}{\begin{tabular}[c]{@{}c@{}}$3.6168$\\ $\cdot 10^{-1}$ \end{tabular} }                                                    & \multicolumn{1}{c|}{\begin{tabular}[c]{@{}c@{}}$-6.7205$\\ $\cdot 10^{-1}$ \end{tabular} }                                                          &   \begin{tabular}[c]{@{}c@{}}$2.2267$\\ $\cdot 10^{-2}$ \end{tabular}                                                         & \multicolumn{1}{c|}{\begin{tabular}[c]{@{}c@{}}$3.7412$\\ $\cdot 10^{-1}$ \end{tabular} }                                                    & \multicolumn{1}{c|}{\begin{tabular}[c]{@{}c@{}}$-6.7573$\\ $\cdot 10^{-1}$ \end{tabular} }                                                        &     \begin{tabular}[c]{@{}c@{}}$2.5913$\\ $\cdot 10^{-2}$ \end{tabular}                                                   \\ \hline
$Z'$ (nd)               & \multicolumn{1}{c|}{\begin{tabular}[c]{@{}c@{}}$1.9606$\\ $\cdot 10^{-1}$ \end{tabular}}                                                    & \multicolumn{1}{c|}{\begin{tabular}[c]{@{}c@{}}$2.0540$\\ $\cdot 10^{-1}$ \end{tabular}}                                                          &   \begin{tabular}[c]{@{}c@{}}$-5.5763$\\ $\cdot 10^{-1}$ \end{tabular}                                                       & \multicolumn{1}{c|}{\begin{tabular}[c]{@{}c@{}}$1.9643$\\ $\cdot 10^{-1}$ \end{tabular}}                                                    & \multicolumn{1}{c|}{\begin{tabular}[c]{@{}c@{}}$1.9973$\\ $\cdot 10^{-1}$ \end{tabular} }                                                        & \begin{tabular}[c]{@{}c@{}}$-5.6094$\\ $\cdot 10^{-1}$ \end{tabular}                                                       \\ \hline
\end{tabular}
\end{table}

The proposed \acrshort{fdc} strategy provides multiple benefits in constellation design and operation leveraging \acrshort{elfos}. While manually adjusting the initial guesses provided by the \acrshort{dam} achieves a desired constellation, the current approach offers a scalable, automatable scheme to reference constellation design. The frequency-domain constraints are seamlessly integrated with other mission requirements, enabling a systematic search for optimal constellations. This strategy also delivers precise control over some frequency coordinates, a task that is difficult to accomplish through manual orbital element adjustments. For example, although the semi-major axis ($\da{a}$ or $a$) primarily influences $\nu_S$, all orbital elements collectively affect the frequency structure. Similarly, achieving a phase separation of $120^\circ$ among three satellites is not naturally achievable solely via orbital element adjustments. Instead, this \acrshort{fdc} technique directly and accurately targets such configurations. Although this application focuses on reference constellation design, the Jacobian matrix $\partial \vec{\mathbb{F}}/\partial \vec{\mathbb{X}}_0$ and its inverse provide valuable operational insights. For example, the impact of stochastic errors in the spacecraft state on downstream phasing information is linearly predicted using the matrix $\partial \vec{\mathbb{F}}/\partial \vec{\mathbb{X}}_0$. This predictive capability enhances operational robustness. 

%% file: text/multiple_shooter.tex
\section{\label{sec:multi}\acrfull{fdc}: Multiple-Shooting Formulation}

Depending on the dynamical regime, the signal $\signalt$ may not be generated from a single, long propagation from an initial state, $\x{0}$. In many environments associated with instability, such a long propagation likely excites nearby hyperbolic manifolds, destroying the quasi-periodic behavior. In such a case, a multiple-shooting formulation typically aids in construction of the long-term trajectories \cite{oguri2020equuleus, bosanac2018trajectory} and retrieval of QPOs. In a similar context, \acrshort{fdc} is extendable to multiple-shooting formulation. Consider the following sequence of state vectors, 
\begin{align}
    \label{eq:x_state} \x{} = \mb \x{0}, \x{1}, \x{2}, ..., \x{\nP-1} \me,
\end{align}
where $\nP$ is the number of patchpoints along the trajectory. The initial epochs (for a non-autonomous case, e.g., the \acrshort{hfem}) corresponding to the patchpoints are predetermined (by assumption) as,
\begin{align}
    \mb \epoch{0}, \epoch{1}, \epoch{2}, ..., \epoch{\nP}   \me,
\end{align}
where the propagation time is also assumed to be fixed as,
\begin{align}
    \Delta \epoch{i_s} = \epoch{i_s+1} - \epoch{i_s}\quad (0 \leq i_s \leq \nP-1).
\end{align}
A typical multiple shooter to ensure the continuity between the patchpoints then enforces,
\begin{align}
    \label{eq:f_state} \f{p} = \mb \xp{0}{\Delta \epoch{0}}{\x{0}}{\epoch{0}} - \x{1} \\ \xp{0}{\Delta \epoch{1}}{\x{1}}{\epoch{1}} - \x{2} \\ ...\\ \xp{0}{\Delta \epoch{\nP-2}}{\x{\nP-2}}{\epoch{\nP-2}} - \x{\nP-1} \me = \vec{0},
\end{align}
where the expression $\xpp$ is defined from Eq. \eqref{eq:xp}. For example, $\xp{0}{\Delta \epoch{0}}{\x{0}}{\epoch{0}}$ denotes the downstream state propagated from the initial state and epoch of $\x{0}$ and $\epoch{0}$, respectively, for the duration of $\Delta \epoch{0}$. Note that Eq. \eqref{eq:f_state} is a vector of length $6(\nP -1)$ while the length of $\x{}$ is $6\nP$. Thus, targeting solely for the continuity between the patchpoints lacks proper boundary constraints, resulting in an underconstrained problem \cite{dei2017trajectory}. The Jacobian matrix for $\f{p}$ with respect to $\x{}$ is constructed as,
\begin{align}
    \label{eq:f_state_ms} \frac{\partial \f{p}}{\partial \x{}} = \mb \stm{0}{\Delta \epoch{0}}{\x{0}}{\epoch{0}} & -\bm{I}_{6\times6} & \bm{0}_{6\times6} & \bm{0}_{6\times6} & ... & \bm{0}_{6\times6} & \bm{0}_{6\times6}  \\
    \bm{0}_{6\times6} &\stm{1}{\Delta \epoch{1}}{\x{1}}{\epoch{1}} & -\bm{I}_{6\times6} & \bm{0}_{6\times6} & ...& \bm{0}_{6\times6} & \bm{0}_{6\times6} \\
    \vdots & \vdots &\vdots &\vdots & \ddots &\vdots & \vdots \\
    \bm{0}_{6\times6} & \bm{0}_{6\times6}& \bm{0}_{6\times6} & \bm{0}_{6\times6}& ... & \stm{0}{\Delta \epoch{\nP-2}}{\x{\nP-2}}{\epoch{\nP-2}} & -\bm{I}_{6\times6} \me.
\end{align}
The extension of \acrshort{fdc} into a multiple-shooting formulation builds on this typical continuity corrector defined with the free variable vector (Eq. \eqref{eq:x_state}), constraint (Eq. \eqref{eq:f_state}), and the Jacobian matrix (Eq. \eqref{eq:f_state_ms}.

\subsection{Sensitivities of the State Vectors on the Frequency Components }

In addition to the continuity constraint $\f{p}$ (Eq. \eqref{eq:f_state}), the frequency-domain constraint $\f{f}$ (Eq. \eqref{eq:frequency_constraint}) is included to formulate \acrshort{fdc} within the multiple-shooting environment. To this end, the sensitivity vectors defined from the single-shooting formulation (Sect. \ref{sec:sensitivity}) are extended to their multiple-shooting formulation counterparts. Specifically, Eq. \eqref{eq:fdt_partial_gomez_required} is now generalized to the multiple-shooting formulation as,
\begin{align}
    \label{eq:fdt_partial_gomez_required_ms} \frac{\partial \freqA{\idxA}}{\partial \x{}}, \frac{\partial \ampA{\idxA}}{\partial \x{}}, \frac{\partial \phaseA{\idxA}}{\partial \x{}},
\end{align}
where a single vector $\x{0}$ from Eq. \eqref{eq:fdt_partial_gomez_required} is extended to the multiple patchpoint state vectors in Eq. \eqref{eq:x_state}. 

\subsubsection{Sensitivities within the \acrfull{laskar} Method}

A similar process from the single-shooting formulation is repeated. Differentiating the constraint vector (Eq. \eqref{eq:laskar_constraint}) and inverting a matrix yields,
\begin{align}
     \label{eq:sensitivity_laskar_ms}  \mb \frac{\partial\freqA{\idxA}}{\partial\x{}} \\ \frac{\partial\ampA{\idxA}}{\partial\x{}} \\ \frac{\partial\phaseA{\idxA}}{\partial\x{}} \me = -\left( \frac{\partial \constraintVL }{\partial \freeVE_\idxA}\right)^{-1}\mb \frac{\partial \constraintC{L, \freqA{}} }{\partial \x{}}\\ \frac{\partial \constraintC{L, \mathcal{C}} }{\partial \x{}}\\ \frac{\partial \constraintC{L, \mathcal{S}} }{\partial \x{}}\me.
\end{align}
where the left side corresponds to a $(3 \times 6\nP)$ matrix. The expressions for the right side are supplied in \hyperref[ap:sensitivities_ms]{Appendix C}.

\subsubsection{Sensitivities within the \acrfull{gomez} Method}

A similar process is now initiated with the \acrshort{gomez} method. Differentiating Eq. \eqref{eq:gomez_constraint_idxA} yields,
\begin{align}
     \label{eq:sensitivity_gomez_ms} \mb \frac{\partial\freqA{\idxA}}{\partial\x{}} \\ \frac{\partial\ampA{\idxA}}{\partial\x{}} \\ \frac{\partial\phaseA{\idxA}}{\partial\x{}} \me = \left(  \frac{\partial \constraintVG }{\partial \freeVE}\right)^{-1}\mb \frac{\partial\cosf(\freqP{\idxA})}{\partial\x{}} \\ \frac{\partial\sinf(\freqP{\idxA})}{\partial\x{}} \\ \frac{d\csf(\freqPO{\idxA})}{\partial\x{}} \me,
\end{align}
deriving the sensitivities for the \acrshort{gomez} approach. Expressions for the right side are supplied in \hyperref[ap:sensitivities_ms]{Appendix C}.

\subsection{Application: Controlling Quasi-Periodic Motion in \acrfull{nrho}}

The capabilities of the multiple shooter formulation within the frequency-domain are demonstrated on a sample scenario that involves the $\rt{9}{2}$ synodic resonant \acrfull{nrho}. This \acrshort{nrho} is envisioned to serve as a baseline orbit for NASA's Gateway mission. Note that a numerical approach already exists in literature \cite{zimovan2023baseline} to supply a long-term baseline solution. However, due to the underconstrained nature of the problem from Eq. \eqref{eq:f_state}, the solution generated within the \acrshort{hfem} may display seemingly ``random'' behavior depending on the numerical process employed. Brown et al. \cite{brown2025hawaii} also investigate the quasi-periodic trajectories that exist in the vicinity of the $\rt{9}{2}$ \acrshort{nrho}, contributing to the expansion of the available solution space. The proposed approach, i.e., \acrshort{fdc}, serves as a complementary approach to pinpoint the solution via explicitly controlling the quasi-periodic motion. 

\begin{figure}[h!]
    \centering
    \begin{subfigure}[b]{0.48\textwidth}
        \centering
        \includegraphics[width = 0.99\textwidth]{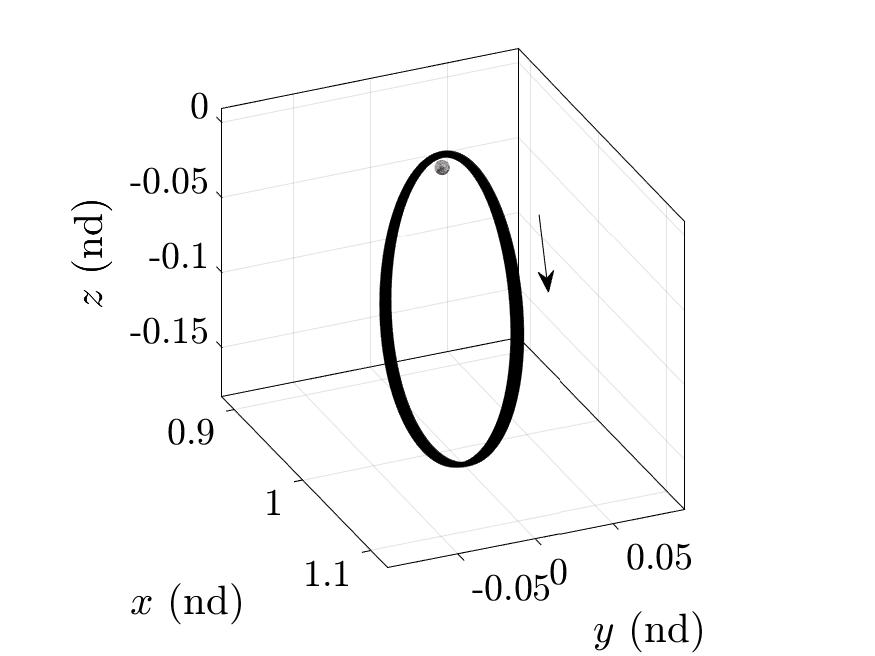}
        \caption{Geometry within the \acrshort{brf}}
        \label{fig:nrho-cr3bp-geometry}
    \end{subfigure}
    \begin{subfigure}[b]{0.48\textwidth}
        \centering
        \includegraphics[width = 0.99\textwidth]{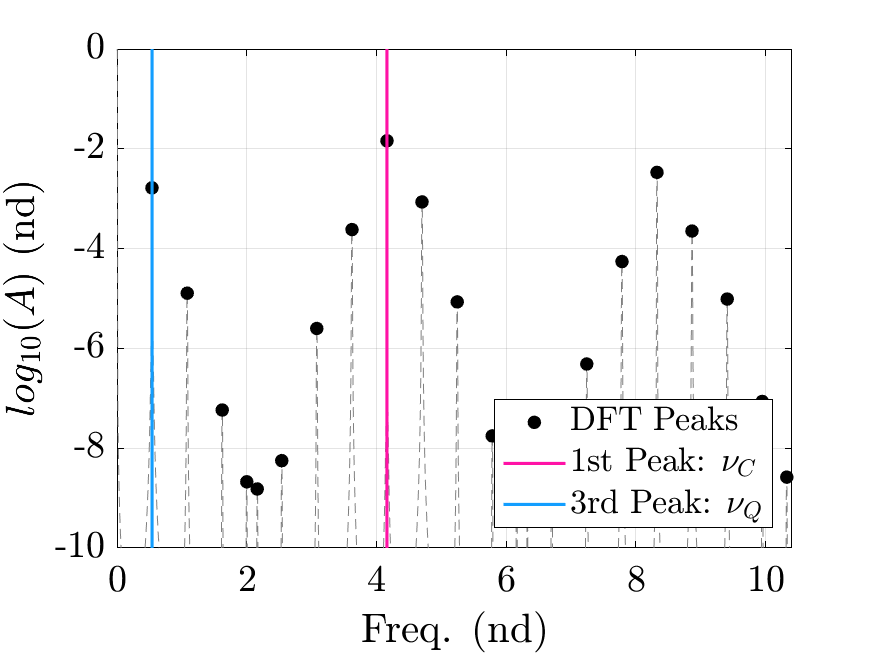}
        \caption{\acrshort{dft} and refined peaks with $\signalt = x$ (\acrshort{laskar})}
        \label{fig:nrho-cr3bp-dft}
    \end{subfigure}
  \caption{Sample \acrshort{cr3bp} quasi-\acrshort{nrho}}
  \label{fig:nrho-cr3bp}
\end{figure}

\begin{figure}[h!]
    \centering
    \begin{subfigure}[b]{0.48\textwidth}
        \centering
        \includegraphics[width = 0.99\textwidth]{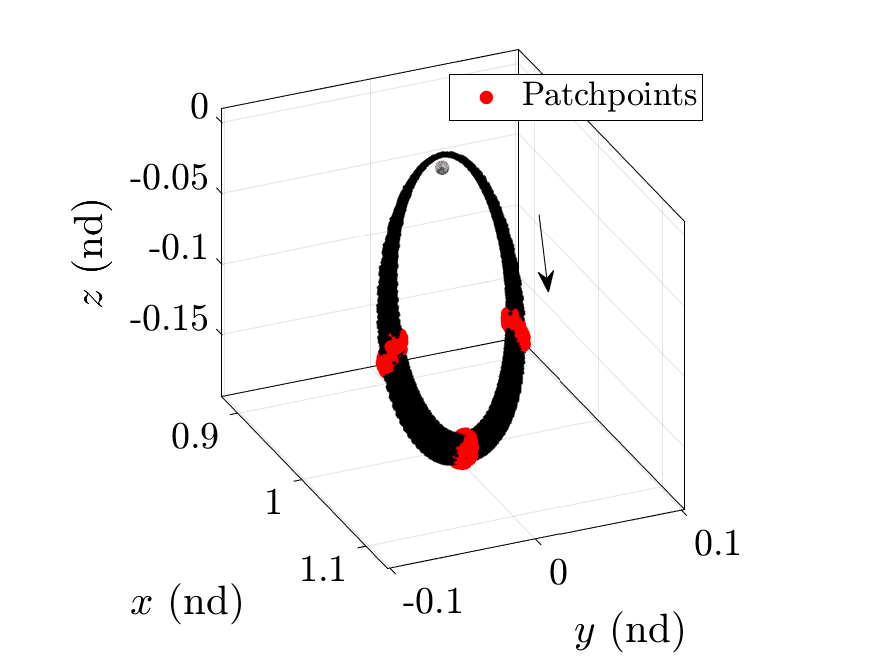}
        \caption{Geometry within the \acrshort{brf}}
        \label{fig:nrho-hfem-geometry}
    \end{subfigure}
    \begin{subfigure}[b]{0.48\textwidth}
        \centering
        \includegraphics[width = 0.99\textwidth]{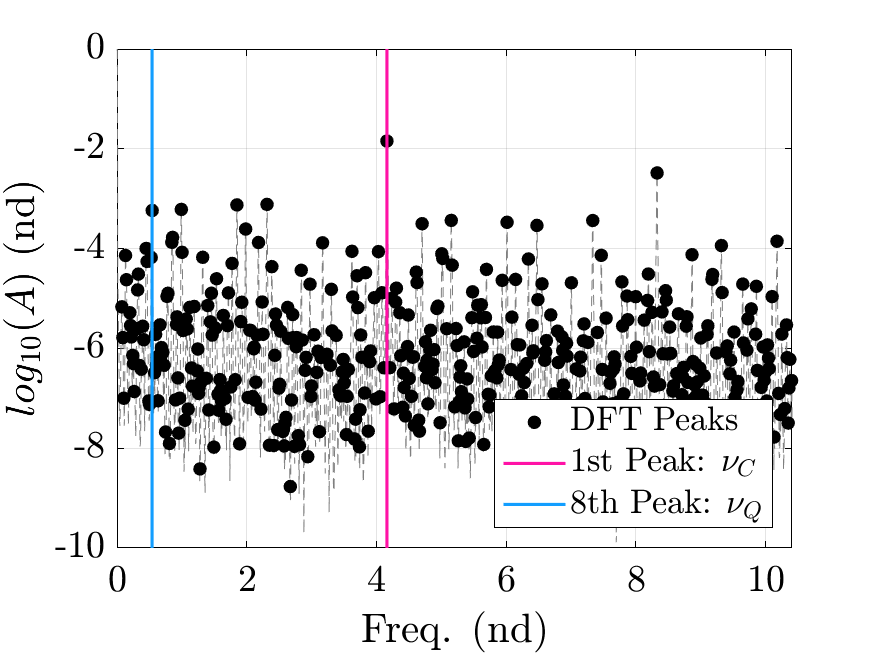}
        \caption{\acrshort{dft} and refined peaks with $\signalt = x$ (\acrshort{laskar})}
        \label{fig:nrho-hfem-dft}
    \end{subfigure}
  \caption{\acrshort{hfem} quasi-\acrshort{nrho} prior to the targeting process}
  \label{fig:nrho-hfem}
\end{figure}

Several dynamical features regarding the $\rt{9}{2}$ \acrshort{nrho} within the \acrshort{cr3bp} are first reviewed, followed by initial observations in the frequency-domain. Within the \acrshort{cr3bp}, the monodromy matrix for the $\rt{9}{2}$ halo admits the following $6$ eigenvalues: $\approx \mb1, 1, -2.1783, -0.4591, 0.6845\pm \imag 0.7290\me$. Note that the unity pair corresponds to the along-orbit and along-family directions. The real reciprocal pair corresponds to the hyperbolic modes, i.e., the stable and unstable modes of the orbit. The complex conjugate pair denotes the oscillatory (center) mode, giving rise to the quasi-periodic motion nearby. Being a $\rt{9}{2}$ \textit{synodic} \acrshort{nrho}, the fundamental frequency corresponding to the underlying period is $\fund{C}=4.1633$ (nd). The second fundamental frequency is $\fund{Q} \approx 0.5415$ (nd), supplied from the center eigenvalues. Due to the existence of the hyperbolic modes, the maximal dimension of the QPO within the \acrshort{cr3bp} for this orbit is $\inter = 2$. Utilizing the GMOS \cite{gomez2001dynamics2, olikara2012numerical} algorithm, a sample 2D-QPO within the \acrshort{cr3bp} is targeted as illustrated in Fig. \ref{fig:nrho-cr3bp-geometry}. The quasi-periodic nature is apparent from the width of the trajectory within the \acrshort{brf}. The corresponding frequency-domain analysis is included in Fig. \ref{fig:nrho-cr3bp-dft} for $\signalt = x$. For this specific case, the first dominant peak corresponds to the $\nu_C$, and the third dominant peak tracks $\nu_Q$. The shape of the \acrshort{nrho} is more visually ``elliptical'' as opposed to the \acrshort{dro} (Fig. \ref{fig:dro-cr3bp}) within the \acrshort{brf}. Thus, more \acrshort{dft} samples $\sampleSize$ aid in preventing aliasing in general as the amplitudes decay at a slower rate for the \acrshort{nrho}. Thus, $\sampleSize = 2^{20}$ is set for the \acrshort{nrho}. Leveraging the QPO from the \acrshort{cr3bp} as an initial guess, a counterpart (continuous) trajectory within the \acrshort{hfem} is supplied as visualized in Fig. \ref{fig:nrho-hfem-geometry}. A multiple-shooting algorithm is employed, similar to the formulation from Eqs. \eqref{eq:x_state} and \eqref{eq:f_state} \cite{zimovan2023baseline}. Over 20 years of total propagation, three patchpoints are placed per revolution as illustrated in Fig. \ref{fig:nrho-hfem-geometry}. The geometry is more dispersed compared to the \acrshort{cr3bp} QPO within \acrshort{brf}. The spectrum also displays complex behavior within the \acrshort{hfem} in Fig. \ref{fig:nrho-hfem-dft}. The fundamental frequencies $\fund{C}$ and $\fund{Q}$ coincide with the 1st and 8th peak, respectively, refined with the \acrshort{laskar} algorithm. The frequency structure for this \acrshort{hfem} quasi-periodic trajectory is illustrated in Table \ref{table:hfem-nrho-freq-initial} with the frequency, amplitude, and phase angle for each fundamental frequency.

\begin{table}[h!]
\centering
\caption{\label{table:hfem-nrho-freq-initial}Initial frequency structure for the \acrshort{nrho} in the \acrshort{hfem} (Fig. \ref{fig:nrho-hfem-dft}}
\begin{tabular}{|l|l|l|}
\hline
 & $\fund{C}$ ($\idxA = 1$)& $\fund{Q}$ ($\idxA = 8$)  \\ \hline
$\freqA{\idxA}$ & $4.163352948064170$ &$0.541347594980843$  \\ \hline
$\ampA{\idxA}$ & $0.014924521456679$ & $0.000614347272201$ \\ \hline
$\phaseA{\idxA}$ & $0.049415700934662$ & $3.280234288629805$ \\ \hline
\end{tabular}
\end{table}

Capability of the proposed \acrshort{fdc} approach in controlling the quasi-periodic motion within the multiple-shooting formulation is demonstrated. An \textit{arbitrary} set of constraints on the frequency structure is formulated as follows,
\begin{align}
    \label{eq:nrho-target-structure} \text{For } \fund{C} : \quad \mb \phaseA{\idxA} - 0 \me = 0, \quad \text{For } \fund{Q} : \quad \mb \ampA{\idxA} - 6.14347272201 \cdot 10^{-5} \\ \phaseA{\idxA} - \pi \me = \vec{0}.
\end{align}
The first type of fundamental frequency, i.e., $\fund{C}$, is extremely challenging to target within the proposed multiple-shooting formulation. As the epochs ($\epoch{}$) are assumed to be fixed, the total propagation time as well as the number of revolutions are fixed; thus, constraining $\fund{C}$ results in an overconstrained problem that displays numerical sensitivities. While a variable-epoch formulation may alleviate these issue, such an avenue remains out-of-scope. For the second type of fundamental frequency, $\fund{Q}$, the amplitude is adjusted to be $1/10$ of the original value as detected in Table \ref{table:hfem-nrho-freq-initial}. The phase angles $\phaseA{\idxA}$ are arbitrarily targeted as $0$ and $\pi$ for $\fund{C}$ and $\fund{Q}$, respectively. These values are selected to be in the vicinity of the initial values (Table \ref{table:hfem-nrho-freq-initial}). As the continuity in the patchpoints is also targeted simultaneously, drastically changing the phase angles often results in numerical sensitivities. Thus, in practice, multiple steps may be required to reach the desired frequency structure. 

\begin{figure}[h!]
    \centering
    \begin{subfigure}[b]{0.48\textwidth}
        \centering
        \includegraphics[width = 0.99\textwidth]{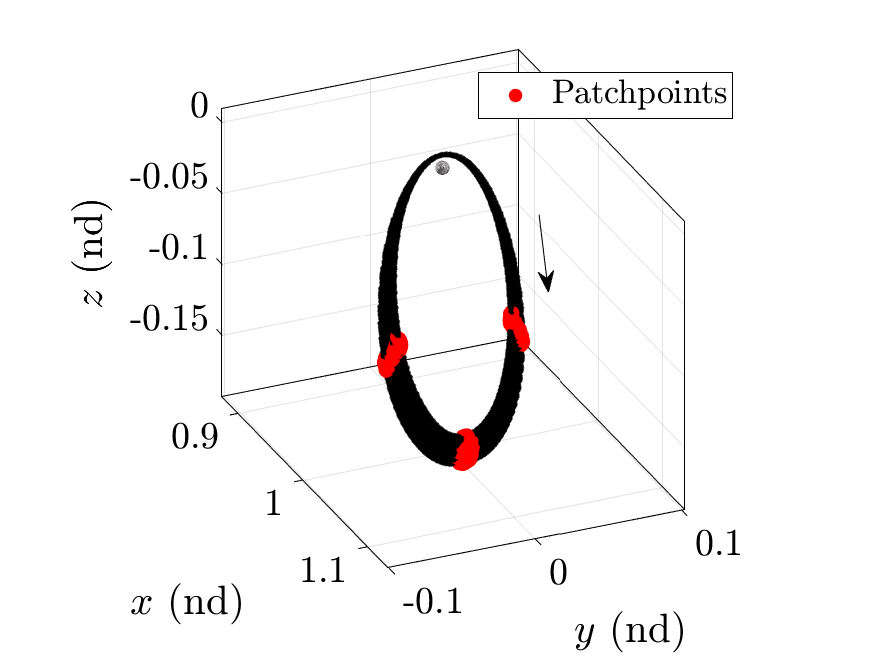}
        \caption{Geometry within the \acrshort{brf}}
        \label{fig:nrho-hfem-converged-geometry}
    \end{subfigure}
    \begin{subfigure}[b]{0.48\textwidth}
        \centering
        \includegraphics[width = 0.99\textwidth]{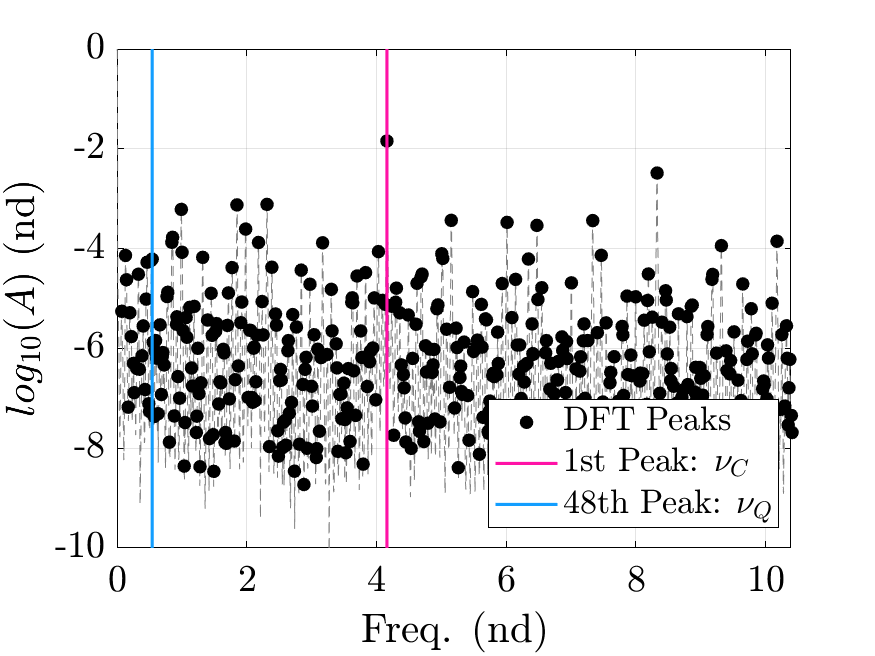}
        \caption{\acrshort{dft} and refined peaks with $\signalt = x$ (\acrshort{laskar})}
        \label{fig:nrho-hfem-converged-dft}
    \end{subfigure}
  \caption{\acrshort{hfem} quasi-\acrshort{nrho} with targeted frequency structure (Eq. \eqref{eq:nrho-target-structure})}
  \label{fig:nrho-hfem-converged}
\end{figure}

Employing the frequency-domain constraints from Eq. \eqref{eq:nrho-target-structure} in addition to the state continuity (Eq. \eqref{eq:f_state}), the initial trajectory from Fig. \ref{fig:nrho-hfem} is adjusted to another solution within the \acrshort{hfem}. The geometry is illustrated in Fig. \ref{fig:nrho-hfem-converged-geometry}. From this view, the visual difference from the initial guess in Fig. \ref{fig:nrho-hfem-geometry} is minimal. However, examining the frequency-domain information from Fig. \ref{fig:nrho-hfem-converged-dft}, it is apparent that the amplitude associated with $\fund{Q}$ is reduced. Now, while $\fund{C}$ is still associated with the most dominant peak, $\fund{Q}$ now tracks the 48th peak within the \acrshort{dft}. The targeted frequency structure is summarized in Table \ref{table:hfem-nrho-freq--target}. As enforced by the constraints from Eq. \eqref{eq:nrho-target-structure}, three frequency components are now successfully targeted. Note that as the amplitude for $\fund{Q}$ decreases, it becomes increasingly challenging to refine it among other peaks as apparent from Fig. \ref{fig:nrho-hfem-converged-dft}. The capabilities of the proposed \acrshort{fdc} for very small $\ampA{\idxA}$ remain an open question. 

\begin{table}[h!]
\centering
\caption{\label{table:hfem-nrho-freq--target}Targeted frequency structure for the \acrshort{nrho} in the \acrshort{hfem} from Fig. \ref{fig:nrho-hfem-converged-dft} (\yes: targeted)}
\begin{tabular}{|l|l|l|}
\hline
 & $\fund{C}$ ($\idxA = 1$)& $\fund{Q}$ ($\idxA = 48$) \\ \hline
$\freqA{\idxA}$ & $4.163383863074772$ &$0.542683258083952$  \\ \hline
$\ampA{\idxA}$ & $0.014924282078229$ & ${0.000061434727220}$ \yes \\ \hline
$\phaseA{\idxA}$ & ${-0.000000000000044}$ \yes& ${3.141592653577986}$ \yes \\ \hline
\end{tabular}
\end{table}

The proposed \acrshort{fdc} approach aids in specifying quasi-periodic trajectories within the \acrshort{hfem}. Functionally, the strategy introduces additional constraints based on the frequency-domain information, reducing the randomness in solution behavior that emanates from the underconstrained nature of the problem \cite{dei2017trajectory}. The approach may be applicable in scenarios where fine-tuning the solution within the \acrshort{hfem} is important. For example, explicitly specifying the phase angle is relevant in avoiding eclipses. Similarly, capabilities to directly control the amplitude and phase angle associated with the quasi-periodic mode ($\fund{Q}$) facilitates the relative trajectory design that leverages QPOs \cite{henry2021expansion, dominguez2024, capannolo2023model}. While the current analysis focuses on the $\rt{9}{2}$ \acrshort{nrho}, the strategy is extendable to other quasi-periodic trajectories associated with inherent instability.

%% file: text/conclusions.tex
\section{\label{sec:summary}Concluding Remarks}

The current study introduces the \acrfull{fdc} strategy for generating multi-dimensional quasi-periodic orbits. Building on insights from existing frequency refinement methods, the approach formulates sensitivities of the frequency structure with respect to the state vector, enabling a frequency-based differential corrections process. The strategy is demonstrated in the Earth-Moon system employing both single- and multiple-shooting formulations. The modularity of the \acrshort{fdc} framework allows flexible integration across a range of dynamical models with evolving fidelity. By incorporating frequency-domain constraints directly into the correction process, the methodology enhances trajectory design capabilities, especially for applications involving constellations or control of oscillatory modes. Future work includes improving numerical robustness and extending the strategy to broader classes of dynamical environments and mission scenarios.

%% file: text/partials.tex
\section*{\label{ap:jacobian}Appendix A: Jacobian Matrix Components for the Refinement Strategies}

\subsection*{\acrshort{laskar}}

Each component of the Jacobian matrix (Eq. \eqref{eq:laskar_jacobian}) for the \acrshort{laskar} is rendered as,
\begin{align}
    \frac{\partial \constraintC{L, \freqAArg}}{\partial \freqA{\idxA}} & = -\frac{1}{2\sqrt{\cosf(\freqA{\idxA})^2 + \sinf(\freqA{\idxA})^2}^3}\left(\cosf(\freqA{\idxA}) \frac{d \cosf(\freqA{\idxA})}{d \freqA{\idxA}} + \sinf(\freqA{\idxA}) \frac{d \sinf(\freqA{\idxA})}{d \freqA{\idxA}}  \right)^2 \nonumber \\
    & + \frac{1}{2\sqrt{\cosf(\freqA{\idxA})^2 + \sinf(\freqA{\idxA})^2}}\left(\left(\frac{d \cosf(\freqA{\idxA})}{d \freqA{\idxA}}\right)^2 + \left(\frac{d \sinf(\freqA{\idxA})}{d \freqA{\idxA}}\right)^2 + \cosf(\freqA{\idxA}) \frac{d^2 \cosf(\freqA{\idxA})}{d \freqA{\idxA}^2} + \sinf(\freqA{\idxA}) \frac{d^2 \sinf(\freqA{\idxA})}{d \freqA{\idxA}^2}  \right)\\
    \frac{\partial \constraintC{L, \freqAArg}}{\partial \ampA{\idxA}} & = 0 \\
    \frac{\partial \constraintC{L, \freqAArg}}{\partial \phaseA{\idxA}} & = 0\\
    \frac{\partial \constraintC{L, \mathcal{C}}}{\partial \freqA{\idxA}} & =  \ampA{\idxA}\cos \phaseA{\idxA} \frac{d \coscosf(\freqA{\idxA})}{d \freqA{\idxA}} - \ampA{\idxA}\sin \phaseA{\idxA}\frac{d \cossinf(\freqA{\idxA})}{d \freqA{\idxA}}  - \frac{d \cosf(\freqA{\idxA})}{d \freqA{\idxA}} \\
    \frac{\partial \constraintC{L, \mathcal{C}}}{\partial \ampA{\idxA}} &= \cos \phaseA{\idxA}\coscosf(\freqA{\idxA})  - \sin \phaseA{\idxA} \cossinf(\freqA{\idxA}) \\
    \frac{\partial \constraintC{L, \mathcal{C}}}{\partial \phaseA{\idxA}} &= -\ampA{\idxA}\sin \phaseA{\idxA}\coscosf(\freqA{\idxA})  - \ampA{\idxA}\cos \phaseA{\idxA} \cossinf(\freqA{\idxA})\\
    \frac{\partial \constraintC{L, \mathcal{S}}}{\partial \freqA{\idxA}} & =  \ampA{\idxA}\cos \phaseA{\idxA}\frac{d \sincosf(\freqA{\idxA})}{d \freqA{\idxA}}  - \ampA{\idxA}\sin \phaseA{\idxA} \frac{d \sinsinf(\freqA{\idxA})}{d \freqA{\idxA}} - \frac{d \sinf(\freqA{\idxA})}{d \freqA{\idxA}} \\ 
    \frac{\partial \constraintC{L, \mathcal{S}}}{\partial \ampA{\idxA}} &= \cos \phaseA{\idxA} \sincosf({\freqA{\idxA}}) - \sin \phaseA{\idxA}\sinsinf({\freqA{\idxA}})   \\
    \frac{\partial \constraintC{L, \mathcal{S}}}{\partial \phaseA{\idxA}} &= -\ampA{\idxA}\sin \phaseA{\idxA} \sincosf{(\freqA{\idxA})} - \ampA{\idxA}\cos \phaseA{\idxA} \sinsinf{(\freqA{\idxA})} .
\end{align}
Note the following additional quantities,
\begin{align}
   & \frac{d^2 \cosf(\freqA{\idxA})}{d \freqA{\idxA}^2} = \frac{2}{\sampleSize}\sum_{\idxT = 0}^{\sampleSize-1} \tidxT^2 \signalt(\tidxT) \han(\idxT) (- \cos (\freqA{\idxA}\tidxT )) \\
   & \frac{d^2 \sinf(\freqA{\idxA})}{d \freqA{\idxA}^2} = \frac{2}{\sampleSize}\sum_{\idxT = 0}^{\sampleSize-1} \tidxT^2 \signalt(\tidxT) \han(\idxT) (- \sin ( \freqA{\idxA}\tidxT)) \\
   & \frac{d \coscosf(\freqA{\idxA})}{d \freqA{\idxA}}  = -\frac{d \sinsinf(\freqA{\idxA})}{d \freqA{\idxA}} =  \frac{2}{\sampleSize}\sum_{\idxT}^{\sampleSize-1}  (-2\tidxT \sin (\freqA{\idxA} \tidxT)   \cos (\tidxT \freqA{\idxA}) \han(\idxT))  \\
   & \frac{d\sincosf(\freqA{\idxA})}{d\freqA{\idxA}} = \frac{2}{\sampleSize}\sum_{\idxT}^{\sampleSize-1}  (-\tidxT (\sin (\freqA{\idxA}) )^2  \han(\idxT)  + \tidxT (\cos (\freqA{\idxA} \tidxT))^2  \han(\idxT)  ).
\end{align}

\subsection*{\acrshort{gomez}}

The components from the Jacobian matrix (Eq. \eqref{eq:gomez_jacobian} are provided as,
\begin{align}
     \label{eq:gomez_partial_1}\frac{\partial \constraintC{G, \mathcal{C}(\freqP{\idxA})}}{\partial \freqA{\idxA}} & = \ampA{\idxA}\cos\phaseA{\idxA} \frac{d \coscosf(\freqP{\idxA})}{d \freqA{\idxA}} - \ampA{\idxA}\sin\phaseA{\idxA} \frac{d \cossinf(\freqP{\idxA})}{d \freqA{\idxA}}\\ 
    \frac{\partial \constraintC{G, \mathcal{C}(\freqP{\idxA})}}{\partial \ampA{\idxA}} & = \cos\phaseA{\idxA} \coscosf(\freqP{\idxA}) - \sin\phaseA{\idxA} \cossinf(\freqP{\idxA}) \\ 
    \frac{\partial \constraintC{G, \mathcal{C}(\freqP{\idxA})}}{\partial \phaseA{\idxA}} & = -\ampA{\idxA}\sin\phaseA{\idxA} \coscosf(\freqP{\idxA}) - \ampA{\idxA}\cos\phaseA{\idxA} \cossinf(\freqP{\idxA}) \\
    \frac{\partial \constraintC{G, \mathcal{S}(\freqP{\idxA})}}{\partial \freqA{\idxA}} & = \ampA{\idxA}\cos\phaseA{\idxA} \frac{d \sincosf(\freqP{\idxA})}{d \freqA{\idxA}} - \ampA{\idxA}\sin\phaseA{\idxA} \frac{d \sinsinf(\freqP{\idxA})}{d \freqA{\idxA}}\\ 
    \frac{\partial \constraintC{G, \mathcal{S}(\freqP{\idxA})}}{\partial \ampA{\idxA}} & = \cos\phaseA{\idxA} \sincosf(\freqP{\idxA}) - \sin\phaseA{\idxA} \sinsinf(\freqP{\idxA}) \\ 
    \frac{\partial \constraintC{G, \mathcal{C}(\freqP{\idxA})}}{\partial \phaseA{\idxA}} & = -\ampA{\idxA}\sin\phaseA{\idxA} \sincosf(\freqP{\idxA}) - \ampA{\idxA}\cos\phaseA{\idxA} \sinsinf(\freqP{\idxA}) \\
    \frac{\partial \constraintC{G, \mathcal{CS}(\freqPO{\idxA})}}{\partial \freqA{\idxA}} & = \ampA{\idxA}\cos\phaseA{\idxA} \frac{d \cscosf(\freqPO{\idxA})}{d \freqA{\idxA}} - \ampA{\idxA}\sin\phaseA{\idxA} \frac{d \cssinf(\freqPO{\idxA})}{d \freqA{\idxA}} \\ 
    \frac{\partial \constraintC{G, \mathcal{CS}(\freqPO{\idxA})}}{\partial \ampA{\idxA}} & = \cos\phaseA{\idxA} \cscosf(\freqPO{\idxA}) - \sin\phaseA{\idxA} \cssinf(\freqPO{\idxA}) \\ 
    \label{eq:gomez_partial_9} \frac{\partial \constraintC{G, \mathcal{CS}(\freqPO{\idxA})}}{\partial \phaseA{\idxA}} & = -\ampA{\idxA}\sin\phaseA{\idxA} \cscosf(\freqPO{\idxA}) - \ampA{\idxA}\cos\phaseA{\idxA} \cssinf(\freqPO{\idxA}).
\end{align}

The ambiguity in $\mathcal{CS}$ for $\freqPO{\idxA}$ is resolved via comparing the norm of the inverse for the Jacobians constructed with $\mathcal{C}$ and $\mathcal{S}$ for $\freqPO{\idxA}$; the one with smaller norm is selected to supply a better conditioned matrix inverse \cite{gomez2010collocation2}.

%% file: text/sensitivities.tex
\section*{\label{ap:sensitivities}Appendix B: Sensitivity Vector Components for the Single-Shooting Formulation}

\subsection*{\acrshort{laskar}}

The vector components in the right-side of Eq. \eqref{eq:fdt_partial_laskar} are supplied as, 
\begin{align}
    \label{eq:laskar_l}\frac{\partial \constraintC{L, \freqAArg}}{\partial \x{0}} & = -\frac{1}{2\sqrt{\cosf(\freqA{\idxA})^2 + \sinf(\freqA{\idxA})^2}^3}L_1L_2 + \frac{1}{2\sqrt{\cosf(\freqA{\idxA})^2 + \sinf(\freqA{\idxA})^2}} \left(L_3+ L_3  \right),
\end{align}
where,
\begin{align}
    L_1 & = \left(\cosf(\freqA{\idxA}) \frac{d \cosf(\freqA{\idxA})}{d \freqA{\idxA}} + \sinf(\freqA{\idxA}) \frac{d \sinf(\freqA{\idxA})}{d \freqA{\idxA}}  \right) \\
    L_2 & = \left(\cosf(\freqA{\idxA}) \frac{\partial \cosf(\freqA{\idxA})}{\partial \x{0}} + \sinf(\freqA{\idxA}) \frac{\partial \sinf(\freqA{\idxA})}{\partial \x{0}}  \right) \\
    L_3 & = \left(\frac{\partial \cosf(\freqA{\idxA})}{\partial \x{0}}\right)\left(\frac{d \cosf(\freqA{\idxA})}{d \freqA{\idxA}}\right) + \left(\frac{\partial \sinf(\freqA{\idxA})}{\partial \x{0}}\right)\left(\frac{d \sinf(\freqA{\idxA})}{d \freqA{\idxA}}\right) \\
    L_4 & = \cosf(\freqA{\idxA}) \frac{\partial}{\partial \x{0}} \left( \frac{d  \cosf(\freqA{\idxA})}{d \freqA{\idxA}}\right) + \sinf(\freqA{\idxA}) \frac{\partial}{\partial \x{0}} \left(\frac{d \sinf(\freqA{\idxA})}{d \freqA{\idxA}}\right).
\end{align}
Additionally,
\begin{align}
    & -\frac{\partial \constraintC{L, \mathcal{C}} }{\partial \x{0}} = \frac{\partial \cosf{(\freqA{\idxA})} }{\partial \x{0}}  \\ 
    & -\frac{\partial \constraintC{L, \mathcal{S}} }{\partial \x{0}} = \frac{\partial \sinf{(\freqA{\idxA})} }{\partial \x{0}},
\end{align}
where,
\begin{align} \frac{\partial  \dft{\signalt}{\freqA{\idxA}}}{\partial  \x{0}} = \frac{1}{N}\sum_{\idxT=0}^{\sampleSize-1} \biggl(  \frac{d \signalt(\tidxT)}{d \x{0}} \han(\idxT)\bigl(\cos( \freqA{\idxA}\tidxT) -\imag  \sin(\freqA{\idxA}\tidxT) \bigr)  \biggr) = \frac{1}{2}\frac{\partial \cosf{(\freqA{\idxA})} }{\partial \x{0}} - \imag  \frac{1}{2} \frac{\partial \sinf{(\freqA{\idxA})} }{\partial \x{0}}. 
\end{align}
Most terms are supplied from Sect. \ref{sec:frequency_refinement} except for,
\begin{align}
     & \frac{\partial}{\partial \x{0}} \left( \frac{d  \cosf(\freqA{\idxA})}{d \freqA{\idxA}}\right) = \frac{2}{\sampleSize}\sum_{\idxT = 0}^{\sampleSize-1} \tidxT \frac{d \signalt(\tidxT)}{d \x{0}} \han(\idxT) (-\sin (\freqA{\idxA}\tidxT))  \\
     & \frac{\partial}{\partial \x{0}} \left(\frac{d \sinf(\freqA{\idxA})}{d \freqA{\idxA}}\right) = \frac{2}{\sampleSize}\sum_{\idxT = 0}^{\sampleSize-1} \tidxT \frac{d \signalt(\tidxT)}{d \x{0}} \han(\idxT) (\cos (\freqA{\idxA}\tidxT)),
\end{align}
with $\frac{d \signalt(\tidxT)}{d \x{0}} = \frac{d \signalt(\tidxT)}{d \xp{0}{\tidxT}{\x{0}}{\epoch{0}}} \stm{0}{\tidxT}{\x{0}}{\epoch{0}}  $, concluding the derivations for the \acrshort{laskar}.  

\subsection*{\acrshort{gomez}}

Evaluation of the right-side of Eq. \eqref{eq:fdt_partial_gomez} requires the following derivative (from Eq. \eqref{eq:dft_signal_bin}),
\begin{align}
    \frac{\partial \dft{\signalt}{\freqP{\idxA}}}{\partial \x{0}} & = \frac{1}{N}\sum_{\idxT=0}^{\sampleSize-1} \biggl(  \frac{d \signalt(\tidxT)}{d \x{0}} \han(\idxT)\bigl(\cos( \freqP{\idxA}\tidxT) -\imag  \sin(\freqP{\idxA}\tidxT) \bigr)  \biggr) \\
    & = \frac{1}{\sampleSize}\sum_{\idxT=0}^{\sampleSize-1} \biggl( \frac{d \signalt(\tidxT)}{d \xp{0}{\tidxT}{\x{0}}{\epoch{0}}} \stm{0}{\tidxT}{\x{0}}{\epoch{0}} \han(\idxT)\bigl(\cos( \freqP{\idxA}\tidxT) -\imag  \sin(\freqP{\idxA}\tidxT) \bigr)  \biggr).
\end{align}
Then, the real and imaginary components constitute the right-side of Eq. \eqref{eq:fdt_partial_gomez}, concluding the derivations.

%% file: text/sensitivities_ms.tex
\section*{\label{ap:sensitivities_ms}Appendix C: Sensitivity Vector Components for the Multiple-Shooting Formulation}

\subsection*{\acrshort{laskar}}

The equation from the single-shooting formultion (Eq. \eqref{eq:laskar_l}) is trivially extended. The process requires following vectors, 
\begin{align}
    \frac{\partial  \dft{\signalt}{\freqA{\idxA}}}{\partial  \x{}} & = \mb\frac{\partial  \dft{\signalt}{\freqA{\idxA}}}{\partial  \x{0}} & \frac{\partial  \dft{\signalt}{\freqA{\idxA}}}{\partial  \x{1}} & ... & \frac{\partial  \dft{\signalt}{\freqA{\idxA}}}{\partial  \x{\nP-1}} \me \\
    \frac{\partial}{\partial \x{}} \left( \frac{d  \cosf(\freqA{\idxA})}{d \freqA{\idxA}}\right) & = \mb \frac{\partial}{\partial \x{0}} \left( \frac{d  \cosf(\freqA{\idxA})}{d \freqA{\idxA}}\right) & \frac{\partial}{\partial \x{1}} \left( \frac{d  \cosf(\freqA{\idxA})}{d \freqA{\idxA}}\right) & ... & \frac{\partial}{\partial \x{\nP}} \left( \frac{d  \cosf(\freqA{\idxA})}{d \freqA{\idxA}}\right)\me \\
    \frac{\partial}{\partial \x{}} \left(\frac{d \sinf(\freqA{\idxA})}{d \freqA{\idxA}}\right) & = \mb \frac{\partial}{\partial \x{0}} \left(\frac{d \sinf(\freqA{\idxA})}{d \freqA{\idxA}}\right) & \frac{\partial}{\partial \x{1}} \left(\frac{d \sinf(\freqA{\idxA})}{d \freqA{\idxA}}\right) & ... & \frac{\partial}{\partial \x{\nP-1}} \left(\frac{d \sinf(\freqA{\idxA})}{d \freqA{\idxA}}\right)\me.
\end{align}
Each element is evaluated as,
\begin{align}
    \label{eq:partial_laskar_ms_element1}\frac{\partial  \dft{\signalt}{\freqA{\idxA}}}{\partial  \x{i_p}} &= \frac{1}{|\mathcal{I}|}\sum_{\idxT \in \mathcal{I} } \biggl(  \frac{d \signalt(\tidxT)}{d\xp{0}{\tidxT-\epoch{i_p}}{\x{i_p}}{\epoch{i_p}}} \stm{0}{\tidxT-\epoch{i_p}}{\x{i_p}}{\epoch{i_p}} \han(\idxT)\bigl(\cos( \freqA{\idxA}\tidxT) -\imag  \sin(\freqA{\idxA}\tidxT) \bigr)  \biggr) \\
     \label{eq:partial_laskar_ms_element2}\frac{\partial}{\partial \x{i_p}} \left( \frac{d  \cosf(\freqA{\idxA})}{d \freqA{\idxA}}\right)& = \frac{2}{|\mathcal{I}|}\sum_{\idxT\in \mathcal{I} }\tidxT  \frac{d \signalt(\tidxT)}{d\xp{0}{\tidxT-\epoch{i_p}}{\x{i_p}}{\epoch{i_p}}} \stm{0}{\tidxT-\epoch{i_p}}{\x{i_p}}{\epoch{i_p}} \han(\idxT) (-\sin (\freqA{\idxA}\tidxT))  \\
     \label{eq:partial_laskar_ms_element3}\frac{\partial}{\partial \x{i_p}} \left(\frac{d \sinf(\freqA{\idxA})}{d \freqA{\idxA}}\right)& = \frac{2}{|\mathcal{I}|}\sum_{\idxT\in \mathcal{I}}\tidxT  \frac{d \signalt(\tidxT)}{d\xp{0}{\tidxT-\epoch{i_p}}{\x{i_p}}{\epoch{i_p}}} \stm{0}{\tidxT-\epoch{i_p}}{\x{i_p}}{\epoch{i_p}} \han(\idxT) (\cos (\freqA{\idxA}\tidxT)),
\end{align}
where the index vector $\mathcal{I}$ is defined as, $ \mathcal{I} = \left\{ i \;\middle|\; \epoch{i_p} \leq \tidxT < \epoch{i_p+1} \right\}$ and $|\mathcal{I}|$ is the length of the vector. Equations \eqref{eq:partial_laskar_ms_element1}-\eqref{eq:partial_laskar_ms_element3} are evaluated for $0 \leq i_p < \nP-1 $. This process concludes the extension of \acrshort{fdc} into the multiple-shooting formulation leveraging the \acrshort{laskar} approach.

\subsection*{\acrshort{gomez}}

Evaluation of the right-side of Eq. \eqref{eq:fdt_partial_gomez} requires the following derivative (from Eq. \eqref{eq:dft_signal_bin}),
\begin{align}
    \frac{\partial \dft{\signalt}{\freqP{\idxA}}}{\partial \x{}} & = \mb \frac{d \dft{\signalt}{\freqP{\idxA}}}{d \x{0}} & \frac{d \dft{\signalt}{\freqP{\idxA}}}{d \x{1}} & ... & \frac{d \dft{\signalt}{\freqP{\idxA}}}{d \x{\nP-1}}\me ,
\end{align}
where,
\begin{align}
    \label{eq:partial_gomez_ms_element}\frac{\partial \dft{\signalt}{\freqP{\idxA}}}{\partial \x{i_p}} & = \frac{1}{|\mathcal{I}|}\sum_{\idxT\in \mathcal{I}} \biggl( \frac{d \signalt(\tidxT)}{d \xp{0}{\tidxT-\epoch{i_p}}{\x{i_p}}{\epoch{i_p}}} \stm{i_p}{\tidxT-\epoch{i_p}}{\x{i_p}}{\epoch{i_p}} \han(\idxT)\bigl(\cos( \freqP{\idxA}\tidxT) -\imag  \sin(\freqP{\idxA}\tidxT) \bigr)  \biggr),
\end{align}
for $0 \leq i_p < \nP-1 $, concluding the derivations.